\baselineskip=14pt plus 2pt
\magnification =1150
\def\sqr#1#2{{\vcenter{\vbox{\hrule height.#2pt\hbox{\vrule width.#2pt
height#1pt\kern#1pt \vrule width.#2pt}\hrule height.#2pt}}}}
\def\square{\mathchoice\sqr64\sqr64\sqr{2.1}3\sqr{1.5}3}

\centerline {\bf Concentration of the invariant measures 
for the periodic}\par
\centerline {\bf Zakharov, KdV, NLS and Gross--Piatevskii equations in $1D$ and $2D$}\par

\vskip.05in
\centerline {Gordon Blower}\par
\centerline {Department of Mathematics and Statistics}\par
\centerline {Lancaster University}\par
\centerline {Lancaster, LA1 4YF}\par
\centerline {England UK}\par
\centerline {email: g.blower$@$ lancaster.ac.uk}\par
\centerline {tel: +44 01524 593962}\par
\centerline {10th April 2015}\par
\vskip.05in
\noindent {\bf Abstract} {This paper concerns  Gibbs measures $\nu$ for some nonlinear PDE over the
$D$-torus ${\bf T}^D$. 
The Hamiltonian $H=\int_{{\bf T}^D} \Vert\nabla u\Vert^2-\int_{{\bf T}^D} \vert u\vert^p$ 
has canonical equations with solutions in
$\Omega_N=\{ u\in L^2({\bf T}^D) :\int \vert u\vert^2\leq N\}$. For $D=1$
and $2\leq p<6$, $\Omega_N$ supports the 
Gibbs measure $\nu (du)=Z^{-1}e^{-H(u)}\prod_{x\in {\bf
T}} du(x)$ which is normalized and formally invariant under the flow generated 
by the PDE. The paper proves that $(\Omega_N, \Vert\cdot\Vert_{L^2}, \nu )$ is a metric probability space of finite
diameter that satisfies the logarithmic Sobolev
inequalities for the periodic $KdV$, the
focussing cubic nonlinear Schr\"odinger equation and the periodic Zakharov system. For suitable subset of 
$\Omega_N$, a logarithmic Sobolev inequality also holds in the critical case 
$p=6$. For
$D=2$, the Gross--Piatevskii equation has 
$H=\int_{{\bf T}^2} \Vert\nabla u\Vert^2-\int_{{\bf T}^2} (V\ast \vert u\vert^2 )
\vert u\vert^2$, for a 
suitable bounded interaction potential $V$
and the Gibbs measure $\nu$ lies on a metric probability space 
$(\Omega , \Vert\cdot\Vert_{H^{-s}}, \nu )$ which satisfies $LSI$. In the above cases, $(\Omega, d, \nu )$ is the limit in
$L^2$ transportation distance of finite-dimensional $(\Omega_n,\Vert\cdot\Vert,\nu_n)$ given by
Fourier sums.}\par 
\par
\noindent {\bf Classification:} 37L55; 35Q53\par
\vskip.05in
\noindent {\bf 1. Introduction}\par
\noindent The periodic Korteweg--de Vries and cubic nonlinear Schr\"odinger equations in space
dimension $D$ may be realised as Hamiltonian systems with an infinite-dimensional phase space
$L^2({\bf T}^D,{\bf R})^{\times 2}$. For instance, the Hamiltonian 
$$H_p(u)={{1}\over {2}}\int_{{\bf T}^D}\Vert \nabla u(\theta )\Vert^2{{d^D\theta}\over{(2\pi
)^D}} -{{\lambda}\over{p}}
\int_{{\bf T}^D}\vert u(\theta )\vert^p
{{d^D\theta}\over{(2\pi )^D}},\eqno(1.1)$$
\noindent is focussing for $\lambda >0$ and defocussing for $\lambda <0$, and the canonical equations
generate the NLS. 
The critical exponent for existence of smooth solutions over all time is $p=2+(4/D)$ by 
[9, p. 6]. In particular $H_4$ generates
the cubic NLS equation for the field $u$. For $N>0$, traditionally called the number operator [15], let $\Omega_N$ be the 
$$\Omega_N=\Bigl\{ u\in L^2({\bf T}^D; {\bf C}):\int_{{\bf T}^D}
\vert u(\theta )\vert^2{{d^D\theta}\over{(2\pi )^D}}\leq
N\Bigr\}.\eqno(1.2)$$
\noindent Observe that $\Omega_N$ is formally invariant under the flow generated by (1.1).\par
\indent For $D=1$, Lebowitz, Rose and Speer [15] introduced an associated Gibbs $\nu$
measure and determined conditions under which $\nu$ can be normalized to define a probability
measure on $\Omega_N$; thus they introduced the modified canonical ensemble as the metric probability
space ${\bf X}=(\Omega_N, \Vert\cdot\Vert_{L^2}, \nu )$. The purpose is to have a
statistical mechanical model of typical solutions of $KdV$ and $NLS$, not just the smooth
solutions. In this paper, we describe concentration of Gibbs measures in terms
of logarithmic Sobolev inequalities, and then use Sturm's theory of metric measure spaces [19] to 
obtain convergence of Gibbs measures on finite-dimensional phase spaces to the
true Gibbs measure.\par
\vskip.05in
\noindent {\bf Definition} ($LSI(\alpha )$) Let $(X,d)$ be a complete and separable metric space, which is
a length space with no isolated points, and $\mu$ a probability measure
on $X$. For $f:X\rightarrow {\bf R}$, introduce the norm of the gradient $\vert\nabla
f(x)\vert=\lim\sup_{y\rightarrow x}
\vert f(y)-f(x)\vert /d(x,y)$. Then $(X, d,\mu )$ satisfies the
logarithmic Sobolev inequality with constant $\alpha >0$ (abbreviated $LSI(\alpha )$) if
$$\int_Xf(x)^2\log \Bigl( f(x)^2/\int_X f^2d\mu \Bigr) \mu (dx)\leq
{{2}\over{\alpha}}\int_X\bigl\vert \nabla f\bigr\vert^2\mu (dx)\eqno(1.3)$$
\noindent for all $f\in L^2(\mu ;X ;{\bf R})$ such that $\vert\nabla
f(x)\vert\in L^2(\mu ;X; {\bf R})$. See [21, chapter 21]. \par
\indent  When $(X, d)=({\bf R}^m, \Vert\cdot
\Vert_{E})$ for some Banach space norm $E$ and $f:{\bf R}^m\rightarrow {\bf R}$ is continuously
differentiable, then we have $\vert \nabla f(x)\vert=\Vert \nabla f(x)\Vert_{E^*}$, 
where $\nabla f$ is the usual gradient and $E^*$ the dual normed space. In the analysis below, we generally apply $LSI(\alpha )$
to functions which may be expressed in terms of the Fourier coordinates, and we require inequalities with constants that
do not depend directly upon the dimension of the phase space. Our results are closely related to those 
of [14], since $LSI$
implies a spectral gap inequality by [21, Theorem 22.28].\par
\indent Bourgain [6] showed that the Gibbs measure on suitably normalized subspaces could 
be constructed from random Fourier
series, so that the Fourier coefficients give an explicit system of canonical coordinates for the
phase space. Let ${\hbox{H}}^s({\bf T}^D)=\{ \sum_{k\in {\bf Z}^D}a_ke^{ik\cdot \theta}: \vert
a_0\vert^2+\sum_{k\in {\bf Z}^d\setminus \{0\}} \vert k\vert^{2s}\vert a_k\vert^2<\infty\}$. Let
$(\gamma_k, \gamma_k')_{k\in {\bf Z}^D}$ be mutually independent standard Gaussian random variables.
Then for $\rho>0$, the periodic Brownian motion
$$b(\theta) =\sum_{k\in {\bf Z}^D} {{(\gamma_k+i\gamma_k')e^{ik\cdot\theta}}\over{\sqrt{\rho +\vert
k\vert^2}}}\qquad (\theta =(\theta_1,\dots ,\theta_D))\eqno(1.4)$$
\noindent lies in ${\hbox{H}}^s({\bf T}^D)$ almost surely for $s<1-(D/2)$.\par 
\indent For $D=1$,
Lebowitz, Rose and Speer [15] showed that for all $N<\infty$ and 
$2\leq p<6$ one can
introduce $Z=Z(N, p,\lambda )>0$ to normalize the Gibbs measure
$$\nu_N(du)=Z^{-1}{\bf I}_{\Omega_N}(u)e^{-H_p(u)}
\prod_{\theta\in {\bf T}} du(\theta )\eqno(1.5)$$
\noindent as a probability on $\Omega_N$. However, for $p>6$, so such $Z$ 
exists. See also [13, 16] for 
alternative constructions of the Gibbs measure.\par
\indent In section 3 of this paper, we prove a
logarithmic Sobolev inequality for $\nu_N$ when $D=1$ and $p=4$. The 
proof depends upon
convexity of the Hamiltonian on $\Omega_N$, and uses a criterion that
 originates with Bakry and Emery [2, 21]. In section 4, we deduce similar results
for the periodic Zakharov system.
In section 5, we use a similar method to prove a $LSI$  for $u\in L^2({\bf T};
{\bf R})$ and $p=3$, where the Hamiltonian generates the KdV equation.
\indent For $D=1$ and $p=6$, there exists $N_0>0$ such that the Gibbs
measure can be normalized on $\Omega_N$ for $N<N_0$, but not for $N>N_0$. In section 6, we obtain a
logarithmic Sobolev inequality for subsets $\Omega_{N,\kappa }=\{u\in \Omega_N: 
\Vert u\Vert_{{\hbox{H}}^s}^2\leq \kappa\}$
and $1/4<s<1/2$ which support most of the Gibbs measure. While
these Gibbs measures are absolutely continuous with respect to Brownian loop,
the Radon--Nikodym derivatives are not logarithmically concave, so our results
do not follow directly from the curvature computations in [19]. Instead we use uniform convexity of
the Hamiltonians on suitable $\Omega_N$, and exploit the property that LSI are stable under suitable perturbations; see [21, Remark 21.5].\par
\indent The partial sums of the spatial Fourier series suggest
classical Hamiltonians on 
finite-dimensional phase spaces $X^n$ given by the low wave numbers, which generate autonomous systems
of ordinary differential equations in the canonical coordinates. Such $X^n$
support Liouville measures $\nu_n$, which are invariant under the flow generated by the
canonical equations, and which give metric probability spaces ${\bf X}^n=(X^n, \Vert
\cdot\Vert_{{\bf R}^{2n}}, \nu_n)$.  We show that for $D=1$ and $p\leq 6$,
the ${\bf X}^n$ converge as metric probability spaces to ${\bf X}$ in the $L^2$ 
transportation distance; this extends the notion of
approximating the solution of a PDE by Fourier partial sums. \par
\indent The lack of smoothness of $b(\theta )$ complicates the analysis 
of the NLS equation in two dimensions, and more drastically in higher space
dimensions. The integral (1.1) with $p=4$ is critical for existence of invariant measures in the 2D focussing case. 
So one introduces a real interaction potential $V$ and works with the Gross--Piatevskii equation
$$i{{\partial u}\over{\partial t}}+{{\partial^2u}\over{\partial\theta_1^2}} 
+{{\partial^2u}\over{\partial \theta_2^2}} +\lambda \bigl( V\ast \vert
u\vert^2\bigr)u=0,\eqno(1.6)$$
\noindent which is also credited to Hartree. In section 7, we impose 
additional hypotheses including $V\in L^\infty({\bf T}^2; {\bf R})$ to obtain
a finite-dimensional logarithmic Sobolev inequality 
and then $V\in {\hbox{H}}^{1+2s}({\bf T}^2; {\bf R})$ to obtain a 
infinite-dimensional LSI. We regard this as realistic, since in their model of a supersolid, Pomeau and Rica [17] consider a soft sphere interaction with $V$
bounded. The Gibbs measure
is supported on distributions in ${\hbox{H}}^{-s}$, so the solutions of (1.6) are typically
not in $L^2({\bf T}^2; {\bf C})$. Nevertheless, in section 8 we achieve convergence in $L^2$ 
transportation distance for finite-dimensional metric probability spaces towards Gibbs measure
on the phase space for the PDE.\par \par
\vskip.05in
\noindent {\bf 2 Metric Measure Spaces for Trigonometric Systems}\par 
\noindent Sturm [19] has developed a theory of metric measure spaces which
refines the metric geometry of Gromov and Hausdorff. We recall some
definitions, which simplify slightly in our setting of probability
spaces, which Sturm calls normalized measure spaces.\par
Let $(X,d)$ be a complete and separable metric space. Now let ${\hbox{Prob}}_0(X)$ be the space of Radon probability measures on $(X,d)$ with the weak
topology; a metric probability space $\hat X$ consists of $(X,d,\mu)$ with 
$\mu\in {\hbox{Prob}}_0(X)$. Suppose that $\mu ,\nu\in {\hbox{Prob}}_0(X)$ and that $\nu$ is absolutely continuous with
respect to $\mu$ and that $f={{d\nu}\over{d\mu}}$ is the Radon--Nikodym derivative. Then the
relative entropy of $\nu$ with respect to $\mu$ is
$${\hbox{Ent}}(\nu\mid\mu )=\int_X f(x)\log f(x)\, \mu (dx),\eqno(2.1)$$
\noindent so that $0\leq {\hbox{Ent}}(\nu\mid\mu )\leq \infty$. For $1\leq s<\infty$, ${\hbox{Prob}}_s(X)$ consists of the subspace
of $\mu\in {\hbox{Prob}}_0(X)$ such that $\int_X \delta (x_0, x)^s\mu (dx)<\infty$ for some or equivalently all
$x_0\in X$. The Wasserstein distance of order $s$ between 
$\mu, \nu\in {\hbox{Prob}}_s(X)$ is
$${\hbox{W}}_s(\mu ,\nu )=\inf_\pi \Bigl\{ \Bigl(\int\!\!\!\int_{X\times X} \delta (x,y)^s\pi (dxdy)\Bigr)^{1/s}: \pi_1=\mu,
\pi_2=\nu\Bigr\}\eqno(2.2)$$
\noindent where $\pi\in {\hbox{Prob}}_s(X\times X)$ with marginals $\pi_1=\mu$ and $\pi_2=\nu $ is called a
transportation plan, and $\delta^s$ is the cost function. Then 
$({\hbox{Prob}}_s(X), {\hbox{W}}_s)$ is a metric space.\par
\indent Suppose further that there exists $\alpha>0$ such that 
$${\hbox{W}}_s(\nu ,\mu )\leq \sqrt{ {{2}\over{\alpha}}{\hbox{Ent}}(\nu\mid\mu )}\eqno(2.3)$$
\noindent for all $\nu\in {\hbox{Prob}}_s(X)$ that are of finite relative entropy with respect to
$\mu$. Then $\mu$ is said to satisfy the transportation inequality $T_s(\alpha )$. We repeatedly
use the result of Otto and Villani that $LSI(\alpha )$ implies 
$T_2(\alpha )$ on Euclidean space; see [21, 22.17].\par
\vskip.05in
\noindent {\bf Definition} ($L^2$ transportation distance) A pseudo metric on a nonempty set $Z$ is a function $\delta
:Z\times Z\rightarrow [0, \infty ]$ that is symmetric, vanishes on the
diagonal, and satisfies the triangle inequality. A coupling of pseudo
metric spaces $(X, \delta_1)$ and $(Y, \delta_2)$ is a pseudo metric space
$(Z, \delta )$ such that $Z=X\sqcup Y$ and $\delta\vert_{X\times X}=\delta_1$
and $\delta\vert_{Y\times Y}=\delta_2$. Given metric probability spaces $\hat X=(X, \delta_1, \mu_1)$ 
and $\hat Y=(Y, \delta_2, \mu_2)$, consider a coupling $\delta$ of these metric spaces and $\pi\in {\hbox{Prob}}_0(X\times Y)$ with marginals $\mu_1$ and
$\mu_2$. Then the $L^2$ transportation distance is
$${\hbox{D}}_{L^2}(\hat X, \hat Y)=\inf_{\delta ,\pi}\Bigl\{
\Bigl(\int\!\!\!\int_{X\times Y}\delta (x,y)^2\pi
(dxdy)\Bigr)^{1/2}\Bigr\},\eqno(2.4)$$
\noindent where the infimum is taken over all such couplings $\delta$ and all
transportation plans $\pi$. One can
easily show that if $\mu_1\in {\hbox{Prob}}_2(X)$ and
$\mu_2\in{\hbox{Prob}}_2(Y),$ then ${\hbox{D}}_{L^2}(\hat X, \hat Y)<\infty $. The diameter
of $\hat X$ is $\sup\{ d(x,y): x,y\in {\hbox{support}}(\mu )\}$. The family of isomorphism classes of metric
probability spaces that have finite diameter gives a metric space
$({\bf X}, {\hbox{D}}_{L^2})$ by results of [19].\par
\indent To obtain $LSI(\alpha )$ for measures on Hilbert space from their finite-dimensional marginals, 
we use the following Lemma, which is related to Theorem 1.3 from [4].\par
\vskip.05in
\noindent {\bf Lemma 2.1} {\sl Let $d\nu=e^{-V(x)}\prod_{j=1}^\infty dx_j$ 
be a Radon probability measure on
$\ell^2({\bf N}; {\bf R})$, and let ${\cal F}_n$ be $\sigma$-algebra that is
generated by the first $n$ coordinate functions, and let $\nu_n$ be the marginal
of $\nu$ for the first $n$ coordinates. Suppose that\par
\indent (i) $V$ is continuously differentiable, and 
$\int \Vert\nabla V(x)\Vert^2_{\ell^2}\nu (dx)<\infty$;\par
\indent (ii) there exists $\alpha>0$ such that $LSI(\alpha )$ holds for 
$X^n=({\bf R}^n, \Vert\cdot\Vert_{\ell^2}, \nu_n)$ for all $n$. \par
\indent Then $LSI(\alpha )$ holds for $X^\infty =
(\ell^2,\Vert\cdot\Vert_{\ell^2},
\nu ),$ and $X^n\rightarrow X^\infty$ in} ${\hbox{D}}_{L^2}$ {\sl as $n\rightarrow\infty.$}\par 
\vskip.05in
\noindent {\bf Proof.} For $0\leq f\in L^2(\ell^2; \nu ; {\bf R})$, let $f_n={\bf E}(f\mid
{\cal F}_n)$, so that $0\leq f_n$ and $f_n\rightarrow f$ almost surely and in $L^2$
as $n\rightarrow\infty$ by the martingale convergence theorem. By Jensen's
inequality applied to the convex function $\varphi (x)=x^2\log x^2$ for $x>0$, we
have 
$$\int f_n^2\log_+f_n^2\, d\nu -\int f_n^2\log_-f_n^2\, d\nu \leq \int f^2\log_+f^2
\,d\nu-
\int f^2\log_-f^2\, d\nu.\eqno(2.5)$$
\noindent Now $\varphi (x)\geq -1/e$, so we can apply the dominated convergence
theorem to the terms with $\log_-$ and Fatou's lemma to the positive terms with
$\log_+$ to deduce that the entropy term on the left-hand side of $LSI$ satisfy 
$$\eqalignno{\int f^2\log\Bigl(f^2/\int f^2d\nu\Bigr)d\nu
&=\lim_{n\rightarrow\infty}
\int f_n^2\log\Bigl(f_n^2/\int f_n^2d\nu\Bigr)d\nu\cr
&\leq\lim\sup_{n\rightarrow\infty} {{2}\over{\alpha}}\int_{{\bf
R}^n}\Vert \nabla f_n(x)\Vert^2_{\ell^2}\nu_n (dx).&(2.6)\cr}$$
\noindent Integrating by parts in the first $n$ coordinates, we see that $\nabla f_n= {\bf E}(\nabla f\mid
{\cal F}_n)+{\bf E}((f_n-f)\nabla V\mid
{\cal F}_n)$, so by the Cauchy--Schwarz inequality
$${{2}\over{\alpha}}\int \Vert\nabla f_n\Vert^2d\nu_n\leq
{{2(1+\varepsilon_n)}\over{\alpha}}\int \Vert\nabla f\Vert^2d\nu
+{{2(1+\varepsilon_n)}\over{\alpha\varepsilon_n }}
\Bigl(\int \vert f_n-f\vert^2d\nu\Bigr)^{1/2}\Bigl(\int \Vert\nabla V\Vert^2
d\nu\Bigr)^{1/2}\eqno(2.7)$$
where we can choose $\varepsilon_n>0$ decreasing to $0$ so that (2.6) and (2.7) give
$$\int f^2\log\Bigl(f^2/\int f^2d\nu\Bigr)d\nu\leq {{2}\over{\alpha}}
\int \Vert\nabla f\Vert^2d\nu.\eqno(2.8)$$
\noindent Hence $\hat X^\infty$ satisfies $LSI(\alpha )$. Now 
$LSI(\alpha )$ implies $T_1(\alpha )$ by [21, 22.17], so\par
\noindent $\int \exp (\alpha\Vert
x\Vert^2/2)\nu (dx)<\infty$. Any continuous and bounded function $f_n:{\bf
R}^n\rightarrow {\bf R}$ may be identified with a function on the first $n$
coordinates of $\ell^2$, so the equation $\int f_nd\nu_n=\int f_nd\nu$ determines
$\nu_n\in {\hbox{Prob}}_2({\bf R}^n)$. We write
$x=(\xi_j)_{j=1}^\infty\in \ell^2$ as 
$x_n=(\xi_1,\dots , \xi_n)$ and $x^n=(\xi_{n+1}, \xi_{n+2},\dots )$ and introduce
$p_n(dx^n\mid \xi_n)\in  {\hbox{Prob}}_2(\ell^2)$ by disintegrating  $\nu (dx)=p_n(dx^n\mid
x_n)\nu_n(dx_n)$ with respect to $\nu_n$; then we couple $X^n$ with $X^\infty$ by mapping $X^n\rightarrow
X^\infty$ via $x_n\mapsto (x_n, 0)$. To transport $\nu_n$ to $\nu$, we select $x_n$
according to the law $\nu_n$, then select $x^n$ according to the law $p_n(dx^n\mid
x_n)$; hence 
$${\hbox{D}}_{L^2}(X^n, X^\infty )^2\leq \int\!\!\!\int_{{\bf R}^n\times \ell^2}\Vert
x^n\Vert^2_{\ell^2}p_n(dx^n\mid x_n)\nu_n(dx_n)=\int_{\ell^2} \Vert x-{\bf E}(x\mid
{\cal F}_n)\Vert^2_{\ell^2}\nu (dx),\eqno(2.9)$$
\noindent which converges to zero as $n\rightarrow\infty$ by the dominated
convergence theorem; so $X^n\rightarrow X^\infty$ in
${\hbox{D}}_{L^2}$ as $n\rightarrow\infty$.\qquad {$\square$}\par

\indent In subsequent sections, we introduce metric probability spaces relating 
to the trigonometric system over ${\bf T}^D$; their properties link curvature, dimension
and the exponent in $H$. In space dimension $D$, let 
$$X^n={\hbox{span}}\{ e^{ik\cdot \theta}: k\in {\bf Z}^D; k=(k_1, \dots ,k_D); \vert
k_j\vert \leq n;\, j=1, \dots ,D \},\eqno(2.10)$$
\noindent so that $\iota_n:X^n\rightarrow X^{n+1}$ is the formal inclusion. When
$n$ is a dyadic power, the metric structure is well described by
Littlewood--Paley theory. For $j\in {\bf N}$, we introduce the dyadic block
$\Delta_j=\{2^{j-1},2^{j-1}+1, \dots , 2^j-1\}$, and for $J=(j_1, \dots ,j_D)\in {\bf N}^D$, let
$\Delta (J) =\Delta_{j_1}\times\dots \times\Delta_{j_D}$. Let $P_J$ be Dirichlet's
projection onto the ${\hbox{span}}\{ e^{ik\cdot\theta}: k\in \Delta (J)\}$, and introduce
the Hamiltonian
$$H_{\Delta (J)}(u)={{1}\over{2}}\int_{{\bf T}^D}\Vert \nabla P_Ju(\theta
)\Vert_{\ell^2}^2{{d^D\theta}\over{(2\pi
)^D}}
-{{\lambda}\over{p}}\int_{{\bf T}^D} \bigl\vert P_J u(\theta )\bigr\vert^p 
{{d^D\theta}\over{(2\pi
)^D}}.\eqno(2.11)$$
\vskip.05in
\noindent {\bf Proposition 2.2} {\sl For $2\leq p\leq 2+(4/D)$ and $N>0$, there exists
$\lambda >0$ such that $H_{\Delta (J)}(u)$ is uniformly convex on
$\Omega_N$.}\par
\vskip.05in
\noindent {\bf Proof.} We observe that $H_{\Delta (J)}$ is twice continuously differentiable on $L^2$,
and
$$\eqalignno{\Bigl({{d^2}\over{dt^2}}&\Bigr)_{t=0} H_{\Delta (J)}(u+tv)\cr
&\geq  
\int_{{\bf T}^D}\Vert \nabla P_Jv(\theta
)\Vert^2_{\ell^2}{{d^D\theta}\over{(2\pi
)^D}}
-\lambda (p-1)\int_{{\bf T}^D} \bigl\vert P_J u(\theta ) \bigr\vert^{p-2}\bigl\vert
P_Jv(\theta )\bigr\vert^2 {{d^D\theta}\over{(2\pi
)^D}}.&(2.12)\cr}$$  
\indent We write $\vert \Delta\vert$ for the cardinality of a finite set $\Delta$, and
observe that by the inequality of the means,
$$\sum_{\ell =1}^D\vert \Delta_{j_\ell}\vert^2\geq D\vert \Delta
(J)\vert^{2/D}.\eqno(2.13)$$
\noindent Hence the first term on the right-hand side of (2.12) satisfies 
$$\int_{{\bf T}^D}\Vert \nabla P_Jv(\theta
)\Vert^2_{\ell^2}{{d^D\theta}\over{(2\pi
)^D}}\geq {{D}\over{4}} \vert \Delta (J)\vert^{2/D}
\int_{{\bf T}^D}\vert P_Jv(\theta )\vert^2{{d^D\theta}\over{(2\pi
)^D}}.\eqno(2.14)$$
\indent Now we introduce de la Vall\'ee Poussin's kernel $K_J$ for $\Delta (J)$, so that 
$\hat K_J(n_1, \dots ,n_D)=1$ for all $(n_1, \dots , n_D)\in \Delta (J)$ and 
$\hat K_J(n_1, \dots ,n_D)=0$ whenever some $n_\ell$ lies outside of
$\Delta_{j_\ell-1}\cup \Delta_{j_\ell}\cup\Delta_{j_\ell+1}.$ Then $P_Ju=K_J\ast P_Ju$,
so by Young's inequality we have constants $c_m$, independent of $u,N$ and
$\Delta_J$ such that 
$$\eqalignno{\int_{{\bf T}^D}\vert P_Ju(\theta )\vert^{2p-4}{{d^D\theta}\over{(2\pi
)^D}}&\leq c_1\bigl\Vert K_J\bigr\Vert^{2p-4}_{L^{(2p-4)/(p-1)}}\Vert
P_Ju\Vert_{L^2}^{2p-4}\cr
&\leq c_2\vert \Delta_J\vert^{p-3}N^{p-2}&(2.15)\cr}$$
\noindent for all $u\in \Omega_N$. Likewise, we have
$$\eqalignno{\int_{{\bf T}^D}\vert P_Jv(\theta )\vert^{4}{{d^D\theta}\over{(2\pi
)^D}}&\leq c_3\bigl\Vert K_J\bigr\Vert^{4}_{L^{4/3}}\Vert P_Ju\Vert^{4}\cr
&\leq c_4\vert \Delta_J\vert \Vert P_Jv\Vert_{L^2}^4.&(2.16)\cr}$$
\noindent Hence by the Cauchy--Schwarz inequality, we have
$$\eqalignno{&\Bigl({{d^2}\over{dt^2}}\Bigr)_{t=0} H_{\Delta (J)}(u+tv)\cr
&\geq  
\Bigl( {{D}\over{4}} \vert \Delta (J)\vert^{2/D}-\lambda (p-1)c_5\vert \Delta
(J)\vert^{(p-2)/2}N^{(p-2)/2}\Bigr)\int_{{\bf T}^D}\vert P_Jv(\theta )\vert^2{{d^D\theta}\over{(2\pi
)^D}},&(2.17)}$$
\noindent where $2/D\geq (p-2)/2$; so given $N>0$, we can choose $\lambda >0$ sufficiently small 
so that the coefficient in parentheses from (2.17) exceeds $D/8$, for all
$J$.\qquad{$\square$}\par
\vskip.05in
\vskip.05in
\noindent {\bf 3. Application to the cubic periodic Schr\"odinger
equation in 1D}\par
\noindent Proposition 2.2 involves an exponent $p=2+(4/D)$ which equals the optimal exponent
for the focussing NLS by [9, page 6]. Such inequalities on dyadic blocks do not of 
themselves lead directly
to $LSI(\alpha )$ on $\Omega_N$. So in sections 3, 4 and 5, we extend Proposition 2.2 to infinite
dimensions. The Hamiltonian
$$H(u)={{1}\over{2}}\int_{{\bf T}}\Bigl\vert
{{\partial u}\over{\partial \theta}}\Bigr\vert^2
{{d\theta}\over{2\pi}}-{{\lambda}\over{4}}\int_{{\bf T}}\bigl\vert
u(\theta )\bigr\vert^4{{d\theta}\over{2\pi}}\eqno(3.1)$$
\noindent may be expressed in terms of the canonical variables $(f,g)$
where $f,g\in L^2([0, 2\pi ]; {\bf R})$, and the field is $u=f+ig$.
Then the canonical equation of motion is the cubic Schr\"odinger equation
$$i{{\partial u}\over{\partial t}}=-{{\partial^2u}\over{\partial
\theta^2}}-\lambda \vert u\vert^2 u,\eqno(3.2)$$
\noindent periodic in $\theta$. Lebowitz, Rose and Speer [15]
considered the Gibbs measures for such partial differential equations,
exploiting the formal invariance of $H(u)$ and the number
operator $N(u)=\int_{{\bf T}}\vert u(\theta )\vert^2 d\theta/(2\pi )$
with respect to time under the flow generated by the NLS. Bourgain [6, 9] introduced a Gibbs
measure $\nu$ for spatially periodic solutions, and established the
existence of a flow for almost all initial data in the support of
$\nu$. \par
\indent Let $(\gamma_j, \gamma_j')_{j=-\infty}^\infty$ be mutually
independent standard Gaussian random variables, so that
$\sum_{j=-\infty; j\neq 0}^\infty e^{ij\theta} (\gamma_j+i\gamma_j')/j$
defines Brownian loop. Let $\lambda , N>0$ and introduce the ball 
$\Omega_N$ as in (1.2). Often it will be more convenient to use the real Fourier
coefficients $a_j,b_j$ of $u$ as canonical coordinates, where $a_j+ib_j=\int
u(\theta ) e^{-ij\theta} d\theta/(2\pi)$. There exists $Z(N, \lambda )>0$ such that 
$$\nu (du)=Z(N, \lambda )^{-1}{\bf
I}_{\Omega_N}(u)\exp\Bigl({{\lambda}\over{4}}\int_{{\bf T}}\bigl\vert
u(\theta )\bigr\vert^4{{d\theta}\over{2\pi}}\Bigr)\prod_{\theta\in [0,
2\pi]} du(\theta ),\eqno(3.3)$$
\noindent defines a probability measure, where as in [15, 6] we define
$$\prod_{\theta\in [0,
2\pi]} du(\theta )=\prod_{j=-\infty; j\neq 0}^{\infty} \exp\Bigl(
-{{j^2}\over{2}}(a_j^2+b^2_j)\Bigr)
{{j^2da_jdb_j}\over{2\pi}},\eqno(3.4)$$  
\noindent namely the measure induced on $L^2$ by Brownian loop. The
indicator ${\bf I}_{\Omega_N}(u)$ restricts the field to the
bounded subset $\Omega_N$ of $L^2$, and ensures convergence. \par
\indent We approximate $\Omega_N$ by finite-dimensional phase spaces. Let
$P_n:L^2\rightarrow {\hbox{span}}\{ e^{ij\theta}: j=-n, \dots ,n\}$ be
the usual Dirichlet projection. Then the Hamiltonian 

$$H_n(u)={{1}\over{2}}\int_{{\bf T}}\Bigl\vert
{{\partial P_nu}\over{\partial \theta}}\Bigr\vert^2
{{d\theta}\over{2\pi}}-{{\lambda}\over{4}}\int_{{\bf T}}\bigl\vert
P_nu(\theta )\bigr\vert^4{{d\theta}\over{2\pi}}\eqno(3.5)$$
\noindent generates the differential equation
$$i{{\partial P_nu}\over{\partial t}}=-{{\partial^2P_nu}\over{\partial
\theta^2}}-\lambda P_n\Bigl( \vert P_nu\vert^2 P_nu\Bigr),\eqno(3.6)$$
\noindent which is associated with a finite-dimensional phase space
$P_nL^2$, and a corresponding Gibbs measure. In terms of the Fourier coefficients, (3.6) is an
autonomous ordinary differential equation. Let $\hat X=(\Omega_N, \Vert\cdot\Vert_{L^2}, \nu
)$ be the metric measure space associated with (3.3), and with $X^n=\Omega_N\cap P_nL^2$, let  
$\hat X^n=(X^n, \Vert\cdot\Vert_{L^2}, \nu_n
)$ be the metric measure space associated with (3.5).\par
\vskip.05in
\noindent {\bf Proposition 3.1} {\sl For $0\leq \lambda N<3/(14\pi^2)$, the Gibbs
measure for $NLS$ on $\Omega_N$ satisfies the logarithmic Sobolev inequality
$$\int_{\Omega_N} F(x)^2 \log\Bigl( F(x)^2/\int  F^2d\nu \Bigr) \nu (dx)\leq
{{2}\over{\alpha}}\int \bigl\Vert \nabla F\bigr\Vert^2_{{\hbox{H}}^{-1}} d\nu,\eqno(3.7)$$
\noindent for $\alpha =1-(14\pi^2N\lambda )/3$.}\par
\vskip.05in 
\noindent {\bf Proof.} For $f=\Re u$ and $g=\Im u$, the Hamiltonian is  
$$H(f+ig)={{1}\over{2}}\int_{\bf T} \Bigl[
\Bigl({{\partial f}\over{\partial \theta}}\Bigr)^2+  
\Bigl({{\partial g}\over{\partial \theta}}\Bigr)^2\Bigr]{{d\theta}\over{2\pi}}-{{\lambda}\over{4}}\int_{\bf T}
\bigl[ f^2+g^2\bigr]^2{{d\theta}\over{2\pi}},\eqno(3.8)$$
\noindent and we aim to show that this is uniformly convex on $\Omega_N$ with
respect to the homogeneous Sobolev norm $(\int \vert f'\vert^2{{d\theta}\over{2\pi}})^{1/2}$ of $\dot
{\hbox{H}}^1$. We
consider $U(f+ig)=\int_{\bf T} (f^2+g^2)^2{{d\theta}\over{2\pi}},$ which contributes a concave term to
the Hamiltonian $H$. We observe that for $0<t<1$ and $f,g,p,q\in {\hbox{H}}^1$,
$$\eqalignno{t&\bigl[f^2+g^2\bigr]^2+(1-t)\bigl[p^2+q^2\bigr]^2-\bigl[
(tf+(1-t)p)^2+(tg+(1-t)q)^2\bigr]^2\cr
&=t(1-t)(f-p)^2\bigl( (1+t+t^2)f^2+(2+2t-2t^2)fp+(2-t+(1-t)^2)p^2\bigr)\cr
&\quad +t(1-t)(g-q)^2\bigl( (1+t+t^2)g^2+(2+2t-2t^2)gq+(2-t+(1-t)^2)q^2\bigr)\cr
&\quad
+2t(1-t)(f-p)(g-q)(f+p)((1+t)g+(1-t)q)\cr
&\quad +2t(1-t)(g-q)^2p^2+2t(1-t)(f-p)^2(tg+(1-t)q)^2.&(3.9)\cr}$$
\noindent We have the basic estimates $\int (f^2+g^2)\leq N$, and likewise $\int
(p^2+q^2)\leq N$, while the Cauchy--Schwarz inequality gives the bounds
$$\Vert f-p\Vert_{L^\infty}^2\leq {{\pi^2}\over{3}}\int\Bigl( {{\partial
f}\over{\partial \theta}}- {{\partial 
p}\over{\partial \theta}}\Bigr)^2 {{d\theta}\over{2\pi}}\eqno(3.10)$$
and likewise for $\Vert g-q\Vert_{L^\infty}$. We integrate (3.9) over ${\bf T}$, and use the ${L^\infty}$ on 
each of the differences $f-p$ and $g-q$ and the squared $L^2$ norm to bound each of the sums; hence
 we have the bound
$$\eqalignno{0&\leq tU(f+ig)+(1-t)U(p+iq)-U(tf+(1-t)p+i(tg+(1-t)q))\cr
&\leq 28Nt(1-t)
\int\Bigl[\Bigl( {{\partial
f}\over{\partial \theta}}- {{\partial 
p}\over{\partial \theta}}\Bigr)^2+\Bigl( {{\partial
g}\over{\partial \theta}}- {{\partial 
q}\over{\partial \theta}}\Bigr)^2\Bigr]  {{d\theta}\over{2\pi}}.&(3.11)\cr}$$
\indent We deduce that $H$ is uniformly convex with respect to the norm 
on $\dot {\hbox{H}}^1$,
with
$$\eqalignno{tH(f+ig)&+(1-t)H(p+iq)-H(tf+(1-t)p+itg+i(1-t)q)\cr
&\geq
t(1-t)\Big({{1}\over{2}}-{{28\lambda N\pi^2}\over{12}}\Bigr)   
\int\Bigl[\Bigl( {{\partial
f}\over{\partial \theta}}- {{\partial 
p}\over{\partial \theta}}\Bigr)^2 +\Bigl( {{\partial
g}\over{\partial \theta}}- {{\partial 
q}\over{\partial \theta}}\Bigr)^2\Bigr]{{d\theta}\over{2\pi}}.&(3.12)\cr}$$
The standard inner product on $L^2({\bf T}^D; d^D\theta/(2\pi )^D; {\bf R})$ is unitarily equivalent to
the standard inner product on $\ell^2({\bf Z}^D)$ under the Fourier transform, and
under this pairing, the dual space of ${\hbox{H}}^s({\bf T}^D; {\bf C})$ is 
${\hbox{H}}^{-s}({\bf T}^D; {\bf C})$. In particular, the dual space of $\dot
{\hbox{H}}^1({\bf T}; {\bf R})$ is $\dot {\hbox{H}}^{-1}({\bf T}; {\bf R})$. So by
Bobkov and Ledoux's Proposition 3.1 of [2], the inequality (3.7) holds for all continuously
differentiable $F:X^n\rightarrow {\bf R}$, which depend on only
finitely many Fourier coefficients. Then by Lemma 2.1, we can deduce
(3.7) for all $F$.\qquad{$\square$}\par
\vskip.05in
\noindent {\bf Theorem 3.2} {\sl Let $p=4$, $D=1$ and $0<N\lambda
<3/(14\pi^2)$. Then $\hat
X^\infty $ of the focussing  cubic
NLS has finite diameter and satisfies $LSI(1-(14\pi^2N\lambda /3))$, and 
$\hat X^n\rightarrow \hat X^\infty$ in} ${\hbox{D}}_{L^2}$ {\sl as $n\rightarrow\infty$.}\par
\vskip.05in
\noindent {\bf Proof.} This follows from Lemma 2.1 and Proposition 3.1. Note
that $\Vert \nabla F\Vert_{{\hbox{H}}^{-1}}\leq \Vert\nabla F\Vert_{L^2}$, so
(3.7) implies (1.3).\qquad{$\square$}\par
\vskip.05in
\noindent {\bf Remark.} One can extend the $L^2$ convergence result in Theorem 3.2
 to all $\lambda ,N>0$,
although the proof becomes more complicated.\par

\vskip.05in
\noindent {\bf 4. Periodic Zakharov system in 1D}\par
\noindent Let $u(\theta,t)$ and $n(\theta,t)$ be periodic in the
space $\theta$ variable; here $u$ is the complex electrostatic envelope field and
$n$ is the real ion density fluctuation field. Then the periodic Zakharov model is the pair of coupled differential
equations 
$$\eqalignno{ i{{\partial u}\over{\partial t}}&=-{{\partial^2u}\over{\partial \theta^2}}+nu;\cr
{{\partial^2 n}\over{\partial t^2}}-{{\partial^2n}\over{\partial
\theta^2}}&={{\partial^2}\over{\partial \theta^2}}\bigl(\vert u\vert^2\bigr).&(4.1)\cr}$$
\noindent The initial condition is 
$$u(\theta,0)=\varphi (\theta),\quad n(\theta,0)=a(\theta),\quad {{\partial n}\over{\partial
t}}(\theta,0)=b(\theta);\eqno(4.2)$$
\noindent and Bourgain [7] established global existence of solutions of (4.1) 
for initial data $\varphi\in {\hbox{H}}^1$,
$a\in L^2$ and $b\in {\hbox{H}}^{-1}.$ We now introduce $V$ as the solution of
$$\eqalignno{{{\partial V}\over{\partial \theta}}&={{\partial n}\over{\partial t}},\cr
{{\partial V}\over{\partial t}}&=-{{\partial n}\over{\partial \theta}}-{{\partial 
}\over{\partial \theta}}\bigl(\vert u\vert^2\bigr),&(4.3)\cr}$$
\noindent such that $\int_{\bf T} V(\theta,t)\, {{d\theta}\over{2\pi}}=0;$ existence may be verified from Fourier series. Then
we introduce the Hamiltonian
$$H(u,n)={{1}\over{4}}\int_{\bf T} \Bigl( 
2\bigl\vert{{\partial u}\over{\partial \theta}}
\bigr\vert^2-\vert u\vert^4+\bigl( n+\vert u\vert^2\bigr)^2+V^2\Bigr)\, {{d\theta}\over{2\pi}},\eqno(4.4)$$
\noindent which suggests that we introduce further variables $\tilde n=(n+\vert
u\vert^2)/\sqrt{2}$ and $W=(d/d\theta )^{-1} V/\sqrt{2}.$ The canonical
variables which lead to the system (4.3) are $(\Re u, \Im u)$ and $(n,
\sqrt{2}W)$. Then $H$ and $\int_{\bf T} \vert u\vert^2\, {{d\theta}\over{2\pi}}$ are 
invariant under the flow, so we can restrict attention to 
$\Omega_B$ as in (1.2) with $D=1$. Then the Gibbs measure on $\Omega_B\times
L^2({\bf T};{\bf R})\times L^2({\bf T};{\bf R})$ is defined by
$$\eqalignno{\nu (dud\tilde ndW)&=Z^{-1}\Bigl[ {\bf I}_{\Omega_B}(u)\exp\Bigl( 
{{1}\over{4}}\int_{\bf T}\vert u\vert^4{{d\theta}\over{2\pi}}
-{{1}\over{2}}\int_{\bf T}\vert{{\partial u}\over{\partial
\theta}}\bigr\vert^2{{d\theta}\over{2\pi}}\Bigr)
\prod_{\theta\in {\bf T}}d^2u(\theta)\Bigr]\cr
&\quad\times \Bigl[ \exp\Bigl( -{{1}\over{2}}\int_{\bf T}  \tilde
n^2\Bigr) \prod_{\theta\in {\bf T}}d\tilde n(\theta )\Bigr] \Bigl[ \exp\Bigl(
-{{1}\over{2}}\int_{\bf T}\bigl( {{\partial
W}\over{\partial\theta}}\bigr)^2{{d\theta}\over{2\pi}}\Bigr)\prod_{\theta\in
{\bf T}}dW(\theta )\Bigr].&(4.5)\cr}$$
\indent We say that $f:L^2\rightarrow {\bf R}$ is a cylindrical function, if there exists a compactly supported
smooth function
$F:{\bf R}^n\rightarrow {\bf R}$ and $\xi_1, \dots ,\xi_n\in L^2$ such that
 $f(\phi )=
F(\langle \phi ,\xi_1\rangle , \dots \langle \phi ,\xi_n\rangle ).$

\vskip.05in
\noindent {\bf Proposition 4.1} {\sl There exists $B>0$ such that the Gibbs measure
for the periodic Zakharov system satisfies
a logarithmic Sobolev inequality for all cylindrical functions.}\par
\vskip.05in
\noindent {\bf Proof.} The Gibbs measure is the direct product of three measures which
satisfy logarithmic Sobolev inequalities, as follows. Let $(\gamma_k)_{k=-\infty}^\infty$ be
mutually independent standard Gaussian random variables, where $\gamma_k$ has distribution $\mu_k$
on ${\bf R}$. Then a typical field $\tilde n$ has the form 
$$\tilde n(\theta ,0)=\sum_{k=1}^\infty \bigl(\gamma_k\cos k\theta
+\gamma_{-k}\sin k\theta \bigr),\eqno(4.6)$$
\noindent which converges in ${\hbox{H}}^{-(1/2)-\varepsilon}$ for all $\varepsilon>0$ almost
surely. By results of Gross and Federbush, each $\mu_k$ satisfies $LSI(1)$ on ${\bf R}$, and likewise $\otimes_{k=-n}^n
\mu_k$ on Eucliean space. The canonical Gaussian measure on $L^2$ has the characteristic property that for
any finite-dimensional subspace $X^n$, the orthogonal projection $P_n:L^2\rightarrow X^n$ induces the
standard Gaussian probability measure on $X^n$ with respect to the induced Euclidean structure; see [18,
page 327]. In particular, this applies to $\otimes_{k=-\infty}^\infty \mu_k$  and the subspace
$X^n={\hbox{span}}\{ \xi_j: j=1, \dots ,n\}$ on which the cylindircal function lives. By [21, page 574; 3] this shows that the middle factor in (4.5) satisfies $LSI(1)$, and there is
no need to truncate the domain  of the $\tilde n$ variable.\par
\indent Likewise, a typical $W$ field initially has the form $W(\theta ,0)=\sum_{k=1}^\infty ( \gamma_k\cos
k\theta +\gamma_{-k}\sin k\theta )/k$ and hence the final factor in (4.5) arises from the direct product of Gaussian
measures that satisfy $LSI(1)$ on ${\bf R}$; hence we have $LSI(1)$ for this product.\par
\indent Finally, the first factor in (4.5) is the Gibbs measure $\nu$ for $NLS$ with $p=4$, so by
Proposition 3.1, $\nu$ satisfies $LSI(1/2)$ for $B<3/28\pi^2$. Combining 
these results,  
as in [21, page 574; 4], we obtain a logarithmic Sobolev inequality where the gradient is
$$\bigl\Vert \nabla F\bigr\Vert^2
=\bigl\Vert \nabla_u F\bigr\Vert_{{\hbox{H}}^{-1}}^2+
\bigl\Vert \nabla_{\tilde n} F\bigr\Vert_{L^2}^2
+\bigl\Vert \nabla_{W} F\bigr\Vert_{{\hbox{H}}^{-1}}^2.\eqno(4.7)$$
\rightline{$\square$}\par
\vskip.05in 
\noindent {\bf 5. Periodic KdV equation in 1D}\par
\noindent Consider $u:{\bf T}\times (0, \infty )\rightarrow {\bf R }$ such that
$u(\,\cdot\, ,t)\in L^2({\bf T})$ for each $t>0$, and introduce the Hamiltonian
$$H(u)={{1}\over{2}}\int_{\bf T} \Bigl( {{\partial u}\over{\partial \theta}}
(\theta,t)\Bigr)^2 {{d\theta}\over{2\pi}} -{{\lambda}\over{6}}\int_{\bf T} u(\theta ,t)^3
{{d\theta}\over{2\pi}},$$
where $\lambda>0$ is the reciprocal temperature. Then the canonical equation of
motion ${{\partial u}\over{\partial t}}=
{{\partial }\over{\partial \theta}}{{\delta H}\over{\delta u}}$ gives the KdV
equation
$${{\partial u}\over{\partial t}}=-{{\partial^3 u}\over{\partial
\theta^3}}-\lambda u{{\partial u}\over{\partial \theta}}.\eqno(5.1)$$
\noindent For a suitably differentiable solution $u$ of (5.1), both $\int
u(\theta,t)^2d\theta/2\pi$ and $H(u)$ are invariant with respect to time. On the ball
$$B_N=\Bigl\{ \phi \in L^2({\bf T}; {\bf R}):\int_{\bf T} \phi (\theta)^2
{{d\theta}\over{2\pi}}\leq N\Bigr\}\eqno(5.2)$$
\noindent with indicator ${\bf I}_{B_N}$ one can define a Gibbs measure
$$\nu (d\phi )=Z_N(\lambda )^{-1}{\bf I}_{B_N}(\phi )e^{-H(\phi )}
\prod_{\theta\in [0,2\pi )} d\phi (\theta)\eqno(5.3)$$
\noindent where $Z_N(\phi )$ is a normalizing constant, chosen to make
$\nu (d\phi )$ a probability measure.\par
\indent The metric probability space $(\Omega_N, \Vert\cdot\Vert_{L^2}, 
\nu )$ arises as
the limit of finite-dimensional metric probability spaces, which are defined in
terms of random Fourier series. Let $X^n=\{ (a_j,b_j)_{j=1}^n \in {\bf R}^{2n}:\phi
(\theta)=\sum_{j=1}^n a_j\cos j\theta +b_j\sin j\theta\in B_N\}$ 
where we introduce the trigonometric polynomial $\phi (\theta)=\sum_{j=1}^n (a_j\cos
j\theta +b_j\sin j\theta)$ and then the probability measure
$$\nu_{n} (dadb)=Z_n^{-1}{\bf I}_{B_N}(\phi )\exp\Bigl(
{{\lambda}\over{6}}\int_{\bf T}\phi
(\theta)^3{{d\theta}\over{2\pi}}\Bigr)\exp\Bigl({-\sum_{j=1}^n j^2(a_j^2+b_j^2)/2}\Bigr)\prod_{j=1}^n da_jdb_j\eqno(5.4)$$
\noindent for a suitable $Z_n=Z_n(N, \lambda )>0$. We then let $\hat X^{(n)}=(X^n,
\Vert\cdot\Vert_{\ell^2},
\nu_{n}),$ which is finite dimensional.\par
\vskip.05in
\noindent {\bf Lemma 5.1} {\sl Suppose that $0\leq \lambda \sqrt{N}<3/\pi^2$. Then the Gibbs
measure satisfies the logarithmic Sobolev inequality}
$$\int_{\Omega_N} f(x)^2\log \Bigl( f(x)^2/\int f^2d\nu\Bigr)\nu (dx)\leq
{{2}\over{\alpha}}\int_{\Omega_N}\bigl\Vert \nabla
f\bigl\Vert^2_{{\hbox{H}}^{-1}}\nu (dx)\eqno(5.5)$$
\noindent {\sl where $\alpha =1-3^{-1}\pi^2\lambda \sqrt{N}$.}\par 
\vskip.05in

\noindent {\bf Proof.} A related result was given in [3] with a larger norm on the right-hand
side. Here we give a proof that is based upon an observation of
Schmuckensl\"ager concerning
uniformly convex Hamiltonians [2, Proposition 3.1]. For $0<t<1$, we have
$$\eqalignno{tH(u)&+(1-t)H(v)-H(tu+(1-t)v)&(5.6)\cr
&={{t(1-t)}\over{2}}\int_{\bf T}\Bigl(
{{\partial u}\over{\partial \theta}} -{{\partial v}\over{\partial
\theta}}\Bigr)^{2}{{d\theta}\over{2\pi}}-{{\lambda t(1-t)}\over{6}}\int_{\bf T}(u-v)^2 \bigl(
(1+t)u+(2-t)v\bigr) {{d\theta}\over{2\pi}}\cr}$$
\noindent where the final term is estimated by the Cauchy--Schwarz
inequality by
$$\eqalignno{\Bigl\vert\int_{\bf T}(u-v)^2 \bigl(
(1+t)u+(2-t)v\bigr) {{d\theta}\over{2\pi}}\Bigr\vert&\leq \Bigl(\int_{\bf
T}(u-v)^4{{d\theta}\over{2\pi}}\Bigr)^{1/2}\Bigl( \int_{\bf T}\bigl
((1+t)u+(2-t)v\bigr)^2{{d\theta}\over{2\pi}}\Bigr)^{1/2}\cr
&\leq  \pi^2\sqrt{N}\int_{\bf T}\Bigl(
{{\partial u}\over{\partial \theta}} -{{\partial v}\over{\partial
\theta}}\Bigr)^{2}{{d\theta}\over{2\pi}}.&(5.7)\cr}$$
\noindent Hence for $\alpha =1-3^{-1}\lambda \pi^2\sqrt{N}>0$, we have a uniformly convex
$H$ such that 
$$tH(u)+(1-t)H(v)-H(tu+(1-t)v)\geq {{t(1-t)\alpha}\over{2}}\int_{\bf T}\Bigl(
{{\partial u}\over{\partial \theta}} -{{\partial v}\over{\partial
\theta}}\Bigr)^{2}{{d\theta}\over{2\pi}};\eqno(5.8)$$
\noindent so $H$ is uniformly convex with respect to
 ${\hbox{H}}^1({\bf T}; {\bf R})$.\qquad{$\square$}\par 
\vskip.05in
\noindent {\bf Theorem 5.2} {\sl Let $0\leq \lambda \sqrt{N}<3/\pi^2$. Then 
 $(\Omega_N, \Vert\cdot\Vert_{L^2}, \nu )$ of $KdV$ has finite diameter, satisfies
$LSI(1-\pi^2\lambda\sqrt{N}/3)$, and is the limit in} ${\hbox{D}}_{L^2}$ {\sl of $\hat
X^n$ as $n\rightarrow\infty.$}\par
\vskip.05in
\noindent {\bf Proof.} Theorem 5.2 follows from lemmas 2.1 and
5.1.\qquad{$\square$}\par
\vskip.05in
\noindent {\bf 6. Logarithmic Sobolev inequality for critical power $p=6$ in 1D}\par
\noindent Now we consider the critical exponent $p=6$, and the Hamiltonian
$$H(u)={{1}\over{2}}\int_{\bf T}\Bigl( {{\partial u}\over{\partial \theta}}\Bigr)^2
{{d\theta}\over{2\pi}}-{{\lambda}\over{6}}\int_{\bf T} u(\theta)^6\, {{d\theta}\over{2\pi}}.\eqno(6.1)$$
Lebowitz, Rose and Speer show that for $0<\lambda \leq 1$, there exists $N_0>0$ such that
the Gibbs measure for $H$ can be normalized on $\Omega_{N}$ for $N<N_0$, but not for $N>N_0$.
To obtain a logarithmic Sobolev inequality, we specialize further and for $1/4<s<1/2$ and
$\kappa{}>0$ let 
$$\Omega_{N,\kappa{}}=\Bigl\{ u\in {\hbox{H}}^s({\bf T}) : \int_{\bf T}\vert u (\theta)\vert^2{{d\theta}\over{2\pi}}\leq N;
\sum_{n=-\infty}^\infty \vert n\vert^{2s}\vert
 \hat u(n)\vert^2\leq \kappa{}\Bigr\}.\eqno(6.2)$$   
\vskip.05in
\noindent {\bf Proposition 6.1} {\sl Let $N<N_0$ and $0<\lambda\leq 1$, and
$1/4<s<1/2$, then let $\nu_N$ be the Gibbs measure on $\Omega_N$ associated with
potential $H$.\par
\indent (i) The sequence of convex and compact subsets
$(\Omega_{N,\kappa{}})_{\kappa{}=1}^\infty$ of 
 $\Omega_N$ is increasing and there exist $\varepsilon, C(\varepsilon )>0$ such that 
$\nu_N(\Omega_{N,\kappa{}})\geq 1-C(\varepsilon )e^{-\varepsilon \kappa{}^2}.$ \par  
\indent (ii) Let $\hat\nu_N$ be $\nu_N$ renormalized on $\Omega_{N, \kappa{}}$
as a probability. Then for all $\kappa{}>0$ there
exists $\alpha =\alpha (\kappa, N)>0$ such that $(\Omega_{N, \kappa{}}, \Vert\cdot\Vert_{L^2},
\hat\nu_N)$ satisfies $LSI(\alpha).$}\par
\vskip.05in 

\noindent {\bf Proof.} (i) Compactness and convexity follow from simple facts about the
Fourier multiplier sequence $(\vert n\vert^{-2s})$ on $L^2$. Let $\mu$ be the 
Gaussian measure on $L^2$ that is induced by Brownian loop. Then by the Cauchy--Schwarz
inequality, we have
$$\int_{\Omega_N} \exp ({\varepsilon \Vert u\Vert^2_{{\hbox{H}}^s}})\nu_N(du)\leq
{{\Bigl(\int_{\Omega_N}\exp ({2\varepsilon \Vert u\Vert^2_{{\hbox{H}}^s}})
\mu (du)
\int_{\Omega_N}\exp ({3^{-1}\lambda \int_{\bf T} u^6})\mu (du)\Bigr)^{1/2}}\over { \int_{\Omega_N} \mu
(du)\int_{\Omega_N} \exp ({6^{-1}\lambda \int_{\bf T}  u^6})\mu (du)}},
\eqno(6.3)$$
\noindent where for suitably small $\varepsilon >0$ the right-hand side integrals are
all finite and together define $C(\varepsilon )$. Then we conclude by applying
Chebyshev's inequality.\par
\indent (ii) For integers $k=1, 2,\dots $, let $\Delta_k=\{ 2^{k-1}, 
2^{k-1}+1,\dots ,
2^{k}-1\}$ be the $k^{th}$ dyadic interval of integers; for $k<0$, let
$\Delta_{k}=\{ n: -n\in \Delta_{-k}\}$; also let $\Delta_0=\{ 0\}$. Next let $K_k$
be de la Vall\'ee Poussin's kernel associated with $\Delta_k$ so $\hat K_k(n)=1$
for all $n\in \Delta_k$, and $\hat K_k(n)=0$ for $n$ outside
$\Delta_{k-1}\cup\Delta_k\cup\Delta_{k+1}$. Also, let
$(\varepsilon_k)_{k=1}^\infty$ be the usual Rademacher functions. By the Littlewood--Paley theorem,
there exist constants $C_1,C_2>0$ {\sl etc.} independent of $u$ such that 
$$\eqalignno{\Vert u\Vert_{L^4}^4&\leq C_1{\bf E}\,
\bigl\Vert\sum_{k=-\infty}^\infty \varepsilon_k K_k\ast u\bigr\Vert^4_{L^4}
\leq C_2\Bigl( \sum_{k=-\infty}^\infty \bigl\Vert
K_k\ast u\bigr\Vert^2_{L^4}\Bigr)^2,&(6.4)\cr}$$
\noindent and we can use Young's inequality to show
$$\bigl\Vert
K_k\ast u\bigr\Vert_{L^4}\leq C_3\Vert K_k\Vert_{L^{4/3}}\Vert K_k\ast
u\Vert_{L^2}\leq C_4\vert \Delta_k\vert^{(1/4)-s}\Vert u\Vert_{{\hbox{H}}^s}.\eqno(6.5)$$
Hence ${\hbox{H}}^s$ embeds continuously in $L^4$. \par
\indent We choose $M>2N_0^2(40^{4s+1}(2\pi \kappa{})^43^{-2})^{1/(4s-1)}$ and introduce
$$U(u)={{M}\over{2}}\int_{\bf T}\vert u(\theta)\vert^2{{d\theta}\over{2\pi}},\eqno(6.6)$$
\noindent so that $U$ is bounded on $\Omega$ with $0\leq U(u)\leq MN\leq MN_0$. Then we
consider the modified Hamiltonian $H(u)+U(u)$, and check that it is uniformly convex, with 
$$\eqalignno{ \Bigl({{d^2}\over{dt^2}}\Bigr)_{t=0}\bigl(
&H(u+tv)+U(u+tv)\bigr)\cr
&=\int_{\bf T}   
\Bigl( {{\partial v}\over{\partial \theta}}\Bigr)^2
{{d\theta}\over{2\pi}}-5\lambda\int_{\bf T}\vert u(\theta)\vert^4\vert v(\theta)\vert^2\, {{d\theta}\over{2\pi}}
+M\int_{\bf T}\vert v(\theta)\vert^2\,{{d\theta}\over{2\pi}}\cr
&\geq \int_{\bf T}   
\Bigl( {{\partial v}\over{\partial \theta}}\Bigr)^2
{{d\theta}\over{2\pi}}-40\lambda \Vert v\Vert_{L^\infty}^2\int_{\bf
T}\Bigl\vert \sum_{k=-\infty}^{-n}K_k\ast u+\sum_{k=n}^\infty
K_k\ast u\Bigr\vert^4 \, {{d\theta}\over{2\pi}}\cr
&\quad +M\int_{\bf T}\vert v(\theta)\vert^2\,
{{d\theta}\over{2\pi}}-40\lambda\Bigl\Vert\sum_{k=-n+1}^{n-1} K_k\ast
u\Bigr\Vert_{L^\infty}^4\int_{\bf T}\vert v(\theta)\vert^2\,
{{d\theta}\over{2\pi}}.&(6.7)\cr}$$
By using the Littlewood--Paley decomposition as above, we obtain the lower bound on (6.7)
$$\bigl( 1-80\lambda \kappa{}^2\vert \Delta_n\vert^{1-4s}\pi^2/3\bigr)\int_{\bf T}   
\Bigl( {{\partial v}\over{\partial \theta}}\Bigr)^2\, {{d\theta}\over{2\pi}}+\bigl( M-40\lambda 
\vert \Delta_n\vert^2
N^2-80\lambda \kappa{}^2\vert \Delta_n\vert^{1-4s})
\int_{\bf T}\vert v(\theta)\vert^2\, {{d\theta}\over{2\pi}}.\eqno(6.8)$$
Now we choose $n$ to be the smallest integer such that 
$2^n=\vert \Delta_n\vert >(160\pi^2\kappa{}^2/3)^{1/(4s-1)}$, so that the first coefficient in
(6.8) exceeds $1/2$, while $M$ was chosen above so that 
$$\Bigl({{d^2}\over{dt^2}}\Bigr)_{t=0}\bigl( H(u+tv)+U(u+tv)\bigr)\geq {{1}\over{2}}
\int_{\bf T}   
\Bigl( {{\partial v}\over{\partial \theta}}\Bigr)^2
{{d\theta}\over{2\pi}}+{{1}\over{2}}\int_{\bf T} v(\theta)^2\, {{d\theta}\over{2\pi}},\eqno(6.9)$$
\noindent and we have uniform convexity. Hence there exists $Z(N)>0$ such that the measure 
$$Z(N)^{-1}e^{-H(u)-U(u)}{\bf I}_{\Omega_{N,\kappa}}(u) 
\prod_{\theta\in [0, 2\pi]} du(\theta)\eqno(6.10)$$
can be normalized and satisfies a logarithmic Sobolev inequality with constant $\alpha_0>0$. The original Gibbs measure 
appears when we perturb the potential by adding the bounded function $U$, to 
remove $-U$; hence by the Holley--Stroock perturbation theorem [11; 21, page 574]
$\nu_N$ also
satisfies a logarithmic Sobolev inequality with constant
$$\alpha \geq \alpha_0 \exp (-NM)\geq\alpha_0 \exp \bigl(
-2(40^{4s+1}(2\pi \kappa{})^43^{-2})^{1/(4s-1)}NN_0^2\bigr).\eqno(6.11)$$ 
\rightline{$\square$}\par
\noindent {\bf 7. The finite-dimensional Gross--Piatevskii equation in 2D}\par
\noindent Let $u\in L^2({\bf T}^2; {\bf C})$, and $a_k+ib_k=\hat u(k)$ be the decomposition of the Fourier coefficients
into real and imaginary parts. With the canonical variables $(a_k, b_k)_{k\in {\bf
Z}^2}$, the Hamiltonian
$$H={{1}\over{2}}\int_{{\bf T}^2}\Vert\nabla
u\Vert^2\,{{d\theta_1}\over{2\pi}}{{d\theta_2}\over{2\pi}}-{{\lambda}\over{4}}\int_{{\bf
T}^2}\bigl( V\ast \vert u\vert^2\bigr) \vert u\vert^2\,
{{d\theta_1}\over{2\pi}}{{d\theta_2}\over{2\pi}}\eqno(7.1)$$
gives rise to the G-P equation (1.5). The $L^2({\bf T}^2; {\bf C})$ norm is invariant
under the flow for smooth periodic solutions.\par
\indent Following Bourgain [8], we introduce a Gibbs measure via random Fourier
series as in (1.4) with $D=2$. Now $b$ does not belong to $L^2({\bf
T}^2;{\bf C})$ almost surely, whereas $b$ defines a distribution in ${\hbox{H}}^{-s}({\bf
T}^2; {\bf C})$ almost surely for all $s>0$. We cannot therefore construct the
canonical ensemble in precisely the same way as in sections 3,4 and 5; instead, we need to
introduce finite-dimensional approximations for which the $L^2$ norms depend upon
the dimension.\par
\indent We define the number operator by 
$$N_n=\sum_{k=(k_1,k_2)\in {\bf Z}^2; \vert k_1\vert, \vert k_2\vert \leq n}{{2}\over{\vert
k\vert^2+\rho}},\eqno(7.2)$$
\noindent so that $N_n\approx 2\log n$ as $n\rightarrow\infty$. Then for $N>0$ let
$\Omega_N$ be as in (1.2) with $D=2$. Let $P_n:L^2({\bf T}^2; {\bf C})\rightarrow {\hbox{span}}\{ e^{ik\cdot\theta};
k\in {\bf Z}^2; k=(k_1, k_2); \vert k_1\vert, \vert k_2\vert \leq n\}$ be 
the usual Dirichlet projection onto the
span of the characters indexed by lattice points in the square of side $2n$ centred at the 
origin. For $B>0$, we let $X^n=P_nL^2\cap
\Omega_{N_n+B}$ with the metric given by the $L^2$ norm, so that the diameter of $X^n$ increases with dimension.
Accordingly, we replace $\vert u\vert^2$ in (7.1) by $\vert u_n\vert^2 -\kappa (N_n+B)$
where $u_n=P_nu$. This is an instance of Wick renormalization.\par
\indent In the following computations, we have integrals over ${\bf T}^2$ with respect to 
$d\theta_1d\theta_2/(2\pi )^2$,
and we suppress the variables of integration. Hence we take the Hamiltonian to be 
$$H_n(u)={{1}\over{2}}\int_{{\bf T}^2} 
\bigl\Vert \nabla u\bigr\Vert^2 -{{\lambda}\over{4}}
\int_{{\bf T}^2} \bigl( V\ast \vert u\vert^2\bigr) \vert u\vert^2 +{{\lambda }\over{2}} \kappa\hat V(0)(N_n+B)\int_{{\bf T}^2}\vert
u\vert^2.\eqno(7.3)$$
We can regard $X^n$ as a compact and convex subset of ${\bf C}^m$ for some $m\leq 4(n+1)^2,$ and define
the Gibbs measure via 
$$\nu_n(dadb)=Z_n^{-1}{\bf I}_{\Omega_{N_n+B}}(u)e^{-H_n(u)}
 \prod_{k=(k_1,k_2)\in {\bf Z}^2; \vert k_1\vert, \vert k_2\vert\leq n}da_kdb_k,\eqno(7.4)$$
\noindent for $u=\sum_{k=(k_1,k_2)\in {\bf Z}^2; \vert k_1\vert, \vert k_2\vert 
\leq n}(a_k+ib_k)e^{ik\cdot\theta}.$ \par
\indent Brydges and Slade [10] consider focussing periodic NLS in 2D and show that some standard routes to renormalization
are blocked. However, allow the
possibility that there exist invariant measures in the case in which $N_n\rightarrow\infty$
and $\lambda_n\rightarrow 0+$ as $n\rightarrow\infty$; see page 489. This is the situation we consider in
Proposition 7.1.\par
\vskip.05in 
\noindent {\bf Proposition 7.1} {\sl (i) Suppose that $V\in L^2({\bf T}^2; {\bf R})$.
 Then for all $B>0$, there exists $\lambda_n>0$ 
such that the Gibbs measure $\nu_n$ on $X^n$ corresponding to $H_n$ satisfies
$LSI(1/2)$, so} 
$$\int_{X^n} f(x)^2\log \Bigl( f(x)^2/\int f^2d\nu_n\Bigr)\nu_n(dx)\leq
4\int_{X^n}\bigl\Vert \nabla
f\bigl\Vert^2_{{\hbox{H}}^{-1}({\bf T}^2)}\nu_n(dx).\eqno(7.5)$$
\indent {\sl (ii) Suppose further that $V\in L^\infty ({\bf T}^2; {\bf R})$ and that $\kappa 
\hat V(0)>3\Vert V\Vert_{L^\infty}$. Then for all $B, 
\lambda >0$ and all $n$, $(X^n,\Vert\cdot\Vert_{L^2},
\nu_n )$ 
satisfies $LSI(1/2)$.}\par
\vskip.05in
\noindent {\bf Proof.}  We prove that the Hamiltonian is uniformly convex, by
introducing
$$\eqalignno{\Bigl( {{d^2}\over{dt^2}}\Bigr)_{t=0} H(u+tw)&
=\int_{{\bf T}^2}\Vert\nabla w\Vert^2 +\lambda \kappa \hat V(0)(N_n+B)
\int_{{\bf T}^2}\vert w\vert^2\cr
&\quad -{{\lambda}\over{2}}\int_{{\bf T}^2}\bigl( \vert w\vert^2 \ast V\bigr)
\vert u\vert^2- {{\lambda}\over{2}}
\int_{{\bf T}^2}\bigl( \vert u\vert^2 \ast V\bigr)\vert w\vert^2\cr
&\quad -{{\lambda}\over{2}}\int_{{\bf T}^2} 
\bigl( (u\bar w+\bar u w)\ast V\bigr) (u\bar w+\bar uw).&(7.6)\cr}$$
\indent (i) By Young's inequality, we have
$$\int_{{\bf T}^2}\bigl( \vert w\vert^2 \ast V\bigr)\vert u\vert^2\leq \bigl\Vert u\bigr\Vert_{L^2}^2\bigl\Vert
V\bigr\Vert_{L^2}\bigl\Vert w\bigr\Vert_{L^4}^2,\eqno(7.7)$$
\noindent and likewise
$$\int_{{\bf T}^2}\bigl( \vert w\vert^2 \ast V\bigr)\vert u\vert^2\leq \bigl\Vert u\bigr\Vert_{L^2}^2\bigl\Vert
V\bigr\Vert_{L^2}\bigl\Vert w\bigr\Vert_{L^4}^2;\eqno(7.8)$$
\noindent while each term in the final term in (7.6) is bounded by Young's inequality and H\"older's inequality, so 
that 
$$\eqalignno{\int_{{\bf T}^2}\bigl( \vert uw\vert \ast \vert V\vert \bigr)\vert uw
\vert&\leq 
\bigl\Vert uw\bigr\Vert_{L^{4/3}}\bigl\Vert
\vert V\vert \ast \vert uw\vert \bigr\Vert_{L^4}\cr
&\leq \bigl\Vert u\bigr\Vert_{L^2}^2\bigl\Vert
V\bigr\Vert_{L^2}\bigl\Vert w\bigr\Vert_{L^4}^2.&(7.9)\cr}$$
\noindent By the Sobolev embedding theorem, we have $\Vert w-\int w\Vert_{L^4}\leq
C_4\Vert \nabla w\Vert_{L^2}$, for some $C_4>0$. Hence
$$\eqalignno{\Bigl( {{d^2}\over{dt^2}}\Bigr)_{t=0} H(u+tw)&
\geq  \Bigl( 1-3\lambda C_4(N_n+B)\Vert
V\Vert_{L^2}\Bigr)\int_{{\bf T}^2}\Vert\nabla w\Vert^2\cr
&\quad +\lambda (N_n+B)\bigl( \kappa \hat V(0)-
3C_4\Vert V\Vert_{L^2}\bigr)\int_{{\bf T}^2} \vert
w\vert^2.&(7.12)\cr}$$ 
\noindent By choosing $\lambda >0$ such that $1/2>3\lambda C_4(N_n+B)\Vert V\Vert_{L^2}$, we
obtain uniform convexity with constant $\alpha =1/2$. Then $LSI(1/2)$
follows from [2, Proposition 3.1].\par 
\indent (ii) When $V$ is bounded, we can use Young's inequality to bound
$$\int \bigl( \vert u\vert^2\ast V\bigr)\vert w\vert^2\leq 
\bigl\Vert V\bigr\Vert_{L^\infty}\bigl\Vert u\bigr\Vert_{L^2}^2 
\bigl\Vert w\bigr\Vert_{L^2}^2,\eqno(7.11)$$
\noindent and likewise for the similar terms in (7.6). Hence we obtain the inequality
$$\eqalignno{\Bigl( {{d^2}\over{dt^2}}\Bigr)_{t=0}& H(u+tw)\cr
&\geq \int_{{\bf T}^2}\Vert\nabla w\Vert^2+
\lambda  \Bigl(\kappa \hat V(0)(N_n+B)-
3\Vert V\Vert_{L^\infty} \int_{{\bf T}^2} \vert u\vert^2\Bigr)
\int_{{\bf T}^2}\vert w\vert^2.&(7.12)\cr}$$
\noindent Again $LSI(\alpha )$ follows from [2, Proposition 3.1].\qquad{$\square$}\par
\vskip.05in
\noindent {\bf 8. The Gross--Piatevskii equation on Sobolev space with negative
index}\par
\noindent To conclude the paper, we obtain a logarithmic Sobolev 
inequality for the G-P
equation (1.5) on a suitable subset of ${\hbox{H}}^{-s}({\bf T}^2; {\bf C}).$ The convolution 
$$\vert u\vert^2\ast V(\theta )=\sum_{m\in {\bf Z}^2}\widehat{ (\vert
u\vert^2)}(m)\hat V(m)e^{im\cdot\theta}\eqno(8.1)$$
\noindent in the potential is to be interpreted probabilistically, since
$u(\theta )=\sum_{k\in {\bf Z}^2\setminus \{
0\}}(\gamma_k+i\gamma_k')e^{ik\cdot
\theta}/\vert k\vert$ does not define an $L^2({\bf T})$ function almost
surely. \par
\indent For $0<s<1/4$, $0<\varepsilon<1/8$, $K_1>0$ and $K_2>5$, let
$$\tilde \Omega =\Bigl\{ (a_j)_{j\in {\bf Z}^2}\in {\bf C}^\infty : 
\sum_{j\in {\bf Z}^2\setminus \{0\}}\vert a_j\vert^2 /\vert
j\vert^{2+2s}\leq K_1;\quad \vert a_j\vert \leq K_2\vert j\vert^{(1/4)-\varepsilon},\quad 
\forall j\in {\bf Z}^2\Bigr\},\eqno(8.2)$$
\noindent so that $\tilde \Omega$ is a convex set. Let $(\gamma_j)_
{j\in {\bf Z}^2}$ be mutually independent standard  complex Gaussian random
variables, so that $\gamma_j$ has distribution $\mu_j$, and let $\tilde \mu$ be
the product measure $\otimes_{j\in {\bf Z}^2}\mu_j$ on ${\bf C}^\infty$. Let
$J:\ell^2({\bf Z}^2; {\bf C})\rightarrow {\hbox{H}}^1({\bf T}^2; {\bf C})$ be the linear map
$J(a_j)=\sum_{j\in {\bf Z}^2\setminus\{0\}}a_je^{ij\cdot\theta}/\vert j\vert$, and let 
$$\Omega = \Bigl\{ u\in {\hbox{H}}^{-s}: \Vert u\Vert_{{\hbox{H}}^{-s}}\leq K_1;\quad  \vert
\hat u(j)\vert\leq K_2\vert j\vert^{-(3/4)-\varepsilon },\quad\forall  j\in {\bf
Z}^2\Bigr\}.\eqno(8.3)$$
\noindent Then $J$ induces a measure $\mu$ on ${\hbox{H}}^{-s}$, which is
mainly supported on $\Omega$.\par
\vskip.05in
\noindent {\bf Theorem 8.1} {\sl Suppose that} $V\in {\hbox{H}}^{1+2s}({\bf
T}^2; {\bf R})$ {\sl for some $s>0$.\par
\indent (i) Then $\mu (\Omega )\rightarrow 1$ as
$K_1, K_2\rightarrow\infty$;\par
\indent (ii) for all $K_1,K_2$ sufficiently large and $0<\varepsilon <1/8$ there
exist $\lambda>0$ and $\alpha>0$ such that the Gibbs measure $\nu$, normalized to
be a probability on $\Omega$,  satisfies $LSI(\alpha)$, so}
$$\int_{\Omega} f(u)^2\log \Bigl( f(u)^2/\int f^2d\nu \Bigr) \nu (du)\leq
{{2}\over{\alpha }}\int_{\Omega} \bigl\Vert \nabla f\bigr\Vert^2_{{\hbox{H}}^{-s}}\, \nu
(du)\eqno(8.4)$$
\noindent {\sl for all $f\in L^2(\Omega ; \nu ; {\bf R})$ that are
differentiable with} $\Vert \nabla f\Vert_{{\hbox{H}}^{-s}}\in L^2(\Omega ; \nu ; {\bf
R})$.\par
\indent {\sl (iii) The transportation cost for cost function}
$c(f,g)=\Vert f-g\Vert^2_{{\hbox{H}}^{-s}}$ {\sl and all} $\omega\in {\hbox{Prob}}_2(\Omega )$
 {\sl that are of finite relative entropy with
respect to $\nu$ satisfies}
$${\hbox{W}}_2(\omega, \nu)^2\leq {{2}\over{\alpha}}{\hbox{Ent}}
(\omega\mid \nu ).\eqno(8.5)$$
\vskip.05in
\noindent {\bf Remark.} The hypotheses imply that $V\in L^\infty$. In
summary, the Gibbs measure produces a metric probability space $(\Omega ,
\Vert\cdot \Vert_{{\hbox{H}}^{-s}}, \nu )$ of finite diameter that
satisfies
$LSI$.\par
\vskip.05in
\noindent {\bf Proof.} (i) We introduce the event
$$\Gamma =\Bigl\{\vert\gamma_j\vert \leq K_2\vert
j\vert^{(1/4)-\varepsilon},
\quad \forall j\in {\bf Z}^2\setminus \{0\}\Bigr\},\eqno(8.6)$$
\noindent which by mutual independence of the $\gamma_j$ has measure 
$$\eqalignno{\tilde\mu (\Gamma )&=\prod_{j\in {\bf Z}^2\setminus \{0\}}\Bigl( 1-2\int_{K_2\vert
j\vert^{1/4-\varepsilon}}^\infty
e^{-s^2/2}{{ds}\over{\sqrt{2\pi}}}\Bigr)\cr
&\geq \exp\Bigl(
-4\sum_{j\in {\bf Z}^2\setminus \{0\}}\int_{K_2\vert
j\vert^{1/4-\varepsilon}}^\infty e^{-s^2/2}{{ds}\over{\sqrt{2\pi}}}\Bigr)
&(8.7)\cr}$$
\noindent since $K_2e^{K_2^2/2}>4$. Also by Chebyshev's inequality, we have
$$\eqalignno{\tilde \mu \Bigl\{\sum_{j\in {\bf Z}^2\setminus \{0\}} {{\vert\gamma\vert^2}\over{\vert
j\vert^{2s+2}}}\geq K^2_1\Bigr\}&\leq e^{-K_1^2/4}\prod_{j\in {\bf Z}^2\setminus \{ 0\}}\Bigl( 1-{{1}\over{2\vert
j\vert^{2+2s}}}\Bigr)^{-1/2}\cr
&\leq \exp \Bigl( -{{K_1^2}\over{4}}+\sum_{j\in {\bf
Z}^2\setminus \{ 0\}} {{1}\over{2\vert j\vert^{2s+2}}}\Bigr);&(8.8)\cr}$$
\noindent so by estimating these sums by the Euler--Maclaurin sum formula, we obtain
$$\tilde \mu (\tilde \Omega )\geq \exp\Bigl(-
{{2(6+\pi )e^{-K_2^2/2}}\over{K_2\sqrt{2\pi}}}\Bigr)-
\exp\Bigl( {{-K_1^2}\over{4}}+{{\pi}\over{2s}}+5\Bigr),\eqno(8.9)$$
\noindent hence $\tilde \mu (\tilde \Omega )\rightarrow 1$ as $K_1,
K_2\rightarrow\infty$.\par 
\indent (ii) By results of Gross and Federbush, each $\mu_j$ satisfies $LSI(1)$
for the standard gradient and distance over ${\bf C}$; hence their
direct product $\tilde \mu$ satisfies $LSI(1)$ on $\tilde \Omega$, where the
norm of the gradient is computed in the norm of $\ell^2$. Lemma 2.1 enables
us to pass from finite to infinite dimensions. We prove below that there exist $\kappa>0$
and $Z>0$ such that $\tilde \nu (da)=Z^{-1}e^{U(J(a))} \tilde\mu (da)$ defines a probability measure
on $\tilde\Omega$ such that 
$$\int_{\tilde \Omega}\exp\Bigl( \kappa\bigl\Vert\nabla (U\circ J)(a)
\bigr\Vert^2_{\ell^2}\Bigr)\tilde\nu
(da)<\infty.\eqno(8.10)$$
\noindent Then $\tilde\nu$ satisfies $LSI(\alpha )$ for some $\alpha>0$ by 
the condition of Aida and Shigekawa [1]; see also [21, Remark 21.5].  Letting $u=J(a)$ and $v=J(b)$, we have
$$\bigl\langle\nabla (U\circ J)(a), b\bigr\rangle_{\ell^2}=(d/dt)_{t=0}U\circ J(a+tb)=\int_{{\bf
T}^2} {{\delta U}\over{\delta u}}(\theta ) v(\theta )\, {{d^2\theta}\over{(2\pi )^2}} ,\eqno(8.11)$$
\noindent while the norms satisfy
$$\eqalignno{\bigl\Vert \nabla (U\circ J)(a)\bigr\Vert_{\ell^2}&=\sup\Bigl\{ \bigl\vert\langle \nabla
(U\circ J)(a),b\rangle\bigr\vert: \Vert b\Vert_{\ell^2}\leq 1\Bigr\}\cr
&= \sup\Bigl\{ \Re \int_{{\bf
T}^2} {{\delta U}\over{\delta u}}(\theta ) v(\theta )\, {{d^2\theta}\over{(2\pi )^2}}: v=J(b); 
\Vert b\Vert_{\ell^2}\leq 1
\Bigr\}\cr
&\leq \Bigl\Vert {{\delta U}\over{\delta u}}\Bigr\Vert_{{\hbox{H}}^{-s}},&(8.12)\cr}$$
\noindent since $J:\ell^2\rightarrow {\hbox{H}}^s$ defines a contractive linear operator for $0<s<1$,
and ${\hbox{H}}^s$ is the dual of ${\hbox{H}}^{-s}$ under the integral pairing.\par
\indent Let $\nu$ be the measure on $\Omega$ that is induced from 
$\tilde \nu$ on $\tilde
\Omega$ by $J$, then normalized to be a probability. Then we obtain the logarithmic Sobolev inequality for the Gibbs measure
$$\eqalignno{ \int_{\Omega} f(\phi )^2\log \Bigl( f(\phi )^2/\int f^2d\nu \Bigr)\nu (d\phi
)&=\int_{\tilde \Omega} f(J(a))^2\log \Bigl( f(J(a))^2/\int f\circ Jd\tilde \nu \Bigr)
e^{U(J(a))}\tilde \mu (da)/Z\cr
&\leq {{2}\over{\alpha}}\int_{\tilde \Omega} \bigl\Vert \nabla (f\circ
J)(a)\bigr\Vert^2_{\ell^2}\tilde \nu (da)\cr
&\leq {{2}\over{\alpha}} \int_{\Omega}\bigl\Vert \nabla f(\phi )\bigr\Vert^2_{{\hbox{H}}^{-s}}\nu
(d\phi),&(8.13)\cr}$$ 
\noindent where the final step follows as in (8.12).\par
\indent So this leaves us with the task of verifying (8.10). The Hamiltonian involves
$$U(u)={{\lambda}\over{4}}\int_{{\bf T}^2}\Bigl( 
\bigl( \vert u\vert^2 -\int
\vert u\vert^2\bigr)\ast V\Bigr)\vert u\vert^2  {{d^2\theta}\over{(2\pi )^2}}\eqno(8.14)$$
\noindent with gradient
$$\eqalignno{ \bigl\langle \nabla U(u),
v\bigr\rangle&=\Bigl({{d}\over{dt}}\Bigr)_{t=0}U(u+tv)&(8.15)\cr
&={{\lambda}\over{4}}\int_{{\bf T}^2}\Bigl[\Bigl( \bigl( \vert u\vert^2 -\int
\vert u\vert^2\bigr)\ast V\Bigr)(u\bar v+v\bar u)+ \Bigl( \bigl( u\bar v+\bar
uv\bigr)\ast V\Bigr)\vert u\vert^2\Bigr] {{d^2\theta}\over{(2\pi )^2}}.\cr}$$
\noindent The integrand involves the Fourier series
$$\bigl( \vert u\vert^2\ast V\bigr)u=\sum_{m\in {\bf Z}^2}\widehat{\bigl( \vert
u\vert^2\bigr)}(m)\hat V(m)\sum_{j\in {\bf Z}^2}\hat u(j)e^{i(j+m)\cdot
\theta},\eqno(8.16)$$
\noindent where $(1+\vert j+m\vert)(1+\vert m\vert)\geq (1+\vert j\vert)$, so
for all $u\in \Omega$ we have
$$\Bigl\Vert \sum_{j\in {\bf Z}^2\setminus\{0\}}\hat
u(j)e^{i(j+m)\cdot\theta}\Bigr\Vert_{{\hbox{H}}^{-s}}\leq \Bigl( \sum_{j\in {\bf
Z}^2\setminus\{0\}}{{\vert
\hat u(j)\vert^2}\over{\vert j\vert^{2s}}}\Bigr)^{1/2}\vert m\vert^s\leq
K_1\vert m\vert^s,\eqno(8.17)$$
\noindent hence
$$\Bigl\Vert \bigl( \vert u\vert^2\ast V\bigr)u\Bigr\Vert_{{\hbox{H}}^{-s}}\leq
K_1\sum_{m\in {\bf Z}^2}\vert m\vert^{s}\vert \hat V(m)\vert \bigl\vert    
\widehat{\bigl( \vert
u\vert^2\bigr)}(m)\bigr\vert .\eqno(8.18)$$
\noindent To estimate the right-hand side of (8.18), we will later use the following lemma.\par
\vskip.05in

\noindent {\bf Lemma 8.2} {\sl (i) The $\widehat{(\vert u\vert^2)}(m)$ are uniformly
exponentially square integrable over $\Omega$ with respect to $\mu$, so there exist $C_1,\kappa>0$ such that}
$$\int_\Omega \exp \Bigl( \kappa^2\bigl\vert \widehat{(\vert u\vert^2)}(m)\bigr\vert^2\Bigr)\mu
(du)<C_1\qquad (m\in {\bf Z}^2\setminus\{ 0\}).\eqno(8.19)$$ 
{\sl (ii) A similar statement holds for $\nu$ on $\Omega$, possibly with different constants.}\par
\vskip.05in
\noindent {\bf Proof.} (i) We have 
$\widehat{(\vert
u\vert^2)}(-m)=\sum_{j}(\gamma_j+i\gamma_j')(\gamma_{j+m}-i\gamma_{j+m}')/\vert
j\vert \vert j+m\vert$, so we require to bound $\sum_{r=1}^\infty d^{(m)}_r$ where each $d^{(m)}_r$ is a
sum over an annulus
$$d^{(m)}_r=\sum_{j\in {\bf Z}^2\setminus \{ 0, -m\}; r-1<\vert j\vert \leq r}
{{\gamma_j\gamma_{j+m}}\over{\vert j\vert\vert j+m\vert}}.\eqno(8.20)$$
\noindent Observe that on $\tilde\Omega$ the random variables $\gamma_j$ are symmetric
and we can independently replace each $\gamma_j$ by $\pm\gamma_j$, without
affecting the distribution of $\tilde\mu$ on $\tilde\Omega$.\par  
\indent The sequence $(d^{(m)}_r)$ is multiplicative in the sense of [12] so that for all
strictly increasing subsequences $r_1<r_2<\dots <r_n$ of integers, 
$$\int_{\tilde\Omega }d^{(m)}_{r_1}d^{(m)}_{r_2}\dots d^{(m)}_{r_n}\,
\tilde\mu (d\gamma )=0.\eqno(8.21)$$
\noindent To see this, consider a product of terms, with one taken from the sum (8.20) for each
factor $d^{(m)}_{r_j}$ and consider the lattice points $\ell$ that index the 
$\gamma_\ell$ from factors in this product. In particular, consider 
 $\ell$ such that the distance from the origin is a maximum, and observe that this is
attained at some point of the form $j+m$, and that $\gamma_{j+m}$ appears only
once in the product, hence integrates to give zero.\par
\indent Observe also that $\vert d^{(m)}_r\vert\leq \delta_r$ where $\delta_r=C_0K_2^2
r^{-(1/2)-\varepsilon}$ for some universal constant $C_0$, so that 
$\delta_r\leq 3C_0^2K_2^4/8\varepsilon$ as follows: 
The most challenging case is when
$\vert m\vert =r$, and we can compare 
$\delta_r\leq K_2^2r^{-(3/4)-\varepsilon}\sum_{j\in {\bf Z}^2\setminus\{ 0,-m\};
r-1\leq \vert j\vert <r}\vert j+m\vert^{-(3/4)-\varepsilon}$ 
with the sum arising with the lattice points $j$ replaced by points equally spaced
around the circle of centre the origin and radius $r$, which produces the integral
 $K_2^2r^{-(1/2)-2\varepsilon}\int_0^{2\pi}\vert \sin (\theta/2)\vert^{-(3/4)-\varepsilon}
d\theta $.\par
\indent Bounded multiplicative
systems satisfy similar concentration inequalities to bounded martingale differences as
 in [20]. 
By Jakubowski and Kwapien's [12] contraction principle, for any convex function $\Phi:{\bf R}^n\rightarrow [0,\infty )$, the inequality 
$$\tilde\mu (\tilde\Omega )^{-1}\int_{\tilde\Omega} \Phi (d^{(m)}_1, \dots ,d^{(m)}_n)\,
d\tilde\mu\leq {\bf E}\Phi
(\delta_1\varepsilon_1, \dots ,\delta_n\varepsilon_n)\eqno(8.22)$$
\noindent holds, where $(\varepsilon_j)_{j=1}^\infty$ is the usual sequence of mutually
independent Rademacher functions. In particular, choosing $\kappa>0$ so that
$\kappa^23C_0^2K_2^4/8\varepsilon<1$, we have
$$\eqalignno{\tilde\mu (\tilde\Omega)^{-1}\int_{\tilde\Omega }\exp\Bigl(
{{\kappa^2}\over{2}} \bigl(\sum_{r=1}^n
d^{(m)}_r\bigr)^2\Bigr) d\tilde \mu &\leq\int_{-\infty}^\infty {\bf
E}\exp\Bigl(t\sum_{r=1}^n
\kappa \delta_r\varepsilon_r\Bigr) \exp (-t^2/2){{dt}\over{\sqrt{2\pi}}}\cr
&=\int_{-\infty}^\infty \prod_{r=1}^n \cosh (\kappa \delta_rt)\, 
\exp (-t^2/2){{dt}\over{\sqrt{2\pi}}}\cr   
&\leq \int_{-\infty}^\infty 
\exp \Bigl({{1}\over{2}}\sum_{r=1}^n \kappa
^2\delta_r^2t^2-{{t^2}\over{2}}\Bigr){{dt}\over{\sqrt{2\pi}}}\cr
&=\Bigl(1-\kappa^2\sum_{r=1}^n \delta_r^2\Bigr)^{-1/2}.&(8.23)\cr}$$
\noindent Letting $n\rightarrow \infty$ and applying Fatou's lemma, we obtain
(8.19).\par
\indent (ii) This follows from (i) by H\"older's inequality.\qquad{$\square$}\par
\vskip.05in 
\noindent {\bf Conclusion of the proof of Theorem 8.1.} (ii) We need to
deduce (8.10) from (8.19). We introduce $C_3>0$ such that
$1/C_3\leq K_1(\pi /s+10)$ such that $C_3\sum_{m\in {\bf
Z}^2\setminus \{0\}}\vert m\vert^{-2-2s}=1$, and then use H\"older's inequality to obtain
$$\eqalignno{\int_\Omega e^{\kappa U(u)}\mu (du)
\leq \prod_{m\in {\bf
Z}^2\setminus\{0\}} &\Bigl[\Bigl( \int_{\Omega}\exp\Bigl[\kappa \bigl\vert \widehat{
(\vert u\vert^2)}(m)\bigr\vert^2/C_3\Bigr]\mu (du)\Bigr)^{C_3\vert \hat
V(m)\vert /2}\cr
& \times\Bigl( \int_{\Omega}\exp\Bigl[\kappa \bigl\vert \widehat{
(\vert u\vert^2)}(-m)\bigr\vert^2/C_3\Bigr]\mu (du)\Bigr)^{C_3\vert \hat
V(m)\vert /2}\Bigr]&(8.24)\cr}$$  
\noindent By Lemma 7.3, all of these integrals converge for sufficiently
small $\kappa >0$, so the Gibbs measure $d\nu =e^Ud\mu$ can be normalized on $\Omega$ to define a
probability measure which is
absolutely continuous with respect to $\mu$.\par 
\indent We can introduce $C(s)\geq  (\pi /s+10)^{-1}$ such that $\sum_{j\in {\bf
Z}^2\setminus\{0\}}C(s)/\vert j\vert^{2+2s}=1$, and then we separate $\hat V$ from
$\widehat {\vert u\vert^2}$ by Cauchy--Schwarz inequality, before applying 
H\"older's inequality to obtain
$$\eqalignno{&\int_\Omega \exp\Bigl[ \kappa_0^2\Bigl(\sum_{m\in {\bf Z}^2\setminus\{0\}}\vert
m\vert^s\vert \hat V(m)\vert \bigl\vert \widehat{(\vert u\vert^2)}(m)\bigr\vert
\Bigr)^2\Bigr]\mu
(du)\cr
&\leq\prod_{m\in {\bf Z}^2\setminus\{0\}}\Bigl( \int_\Omega \exp\Bigl[
{{\kappa_0^2}\over{C(s)}}\sum_{j\in {\bf
Z}^2\setminus\{0\}}\vert j\vert^{2+4s}\vert \hat V(j)\vert^2 \bigl\vert
\,\,\widehat{(\vert
u\vert^2)}(m)\bigr\vert^2\Bigr]\mu (du)\Bigr)^{C(s)/\vert
m\vert^{2+2s}}.&(8.25)\cr}$$
\noindent By taking $\kappa_0>0$ sufficiently small, we can ensure that all the
integrals and the product converge. This confirms that (8.10) holds, and hence gives the
logarithmic Sobolev inequality.\par
\indent (iii) The transportation inequality follows from the logarithmic Sobolev
inequality (8.4) as in [21, Theorem 22.17].\qquad{$\square$}\par
\vskip.05in
\indent Let $X^n={\hbox{span}}\{ e^{ij\cdot\theta}: j\in {\bf Z}^2; \vert j\vert\leq n\}$ be
the subspace of $L^2({\bf T}^2; {\bf C})$ that is spanned by the characters
that are indexed by the lattice points in the disc with radius $n$, and let 
$P_n:L^2({\bf T}^2; {\bf C})\rightarrow X^n$ be the orthogonal projection.
Let $\nu_n$ be the Gibbs measure
$$\nu_n(du)=Z_n^{-1}{\bf
I}_{\Omega}(u)\exp\bigl(U(P_nu)\bigr)\prod_{m\in {\bf
Z}^2\setminus\{ 0\}}e^{-\vert
m\vert^2(a^2_m+b^2_m)/2}da_mdb_m/2\pi.\eqno(8.26)$$
\noindent Let $\omega_n$ be the marginal distribution of $\nu_n$ on $X^n$.\par
\vskip.05in
\noindent {\bf Corollary 8.3} {\sl The } 
$(X^n\cap \Omega , \Vert\cdot\Vert_{{\hbox{H}}^{-s}}, \omega_n)$ {\sl converge in}
${\hbox{D}}_{L^2}$ {\sl to} 
$(\Omega , \Vert\cdot\Vert_{{\hbox{H}}^{-s}}, \nu )$ {\sl as
$n\rightarrow\infty$.}\par
\vskip.05in
\noindent {\bf Proof.} (i) First we prove that $U(P_nu)\rightarrow U(u)$ almost surely and in $L^2$ with respect to $\mu$ on $\Omega$ as
$n\rightarrow\infty$. The difference in the potentials has a Fourier expansion
$$\eqalignno{U(P_nu)-U(u)&=\int (V\ast \vert P_nu\vert^2)\vert P_nu\vert^2 -
\int (V\ast \vert u\vert^2)\vert u\vert^2\cr
&=\sum_{m}\hat V(m)\Bigl( \widehat{\bigl(\vert
P_nu\vert^2\bigr)}(m)-\widehat{\bigl(\vert
u\vert^2\bigr)}(m)\Bigr)\Bigl(\widehat{\bigl( \vert u\vert^2\bigr)}(-m)\Bigr)\cr
&\quad + \sum_{m}\hat V(m)\Bigl(\widehat{\bigl( \vert P_nu\vert^2\bigr)}(m)\Bigr)
\Bigl(\widehat{\bigl(\vert
P_nu\vert^2\bigr)}(-m)- \widehat{\bigl(\vert u\vert^2\bigr)}(-m)\Bigr);&(8.27)\cr}$$
\noindent hence
$$\eqalignno{\bigl\vert &U(P_pu)-U(P_nu)\bigr\vert&(8.28)\cr
&\leq 2\sum_{m\in {\bf Z}^2} \vert \hat V(m)\vert
\bigl\vert \widehat{(\vert P_pu\vert^2)}(m)- \widehat{(\vert P_nu\vert^2)}(m)\bigr\vert
\Bigl(  \bigl\vert \widehat{(\vert P_pu\vert^2)}(m)- \widehat{(\vert P_nu\vert^2)}(m)\bigr\vert
+\bigl\vert \widehat{(\vert u\vert^2)}(m)\bigr\vert\Bigr)\cr}$$
\noindent where 
$$\bigl\vert \widehat{(\vert P_pu\vert^2)}(m)- \widehat{(\vert P_nu\vert^2)}(m)\bigr\vert\leq 
\Bigl\vert \sum_{r=n+1}^\ell d_r^{(m)}\Bigr\vert .\eqno(8.29)$$   
\noindent We observe that $\widehat{(\vert P_nu\vert^2)}(m)$ has a similar
expansion to (8.31), except that only those $j$ with $\vert j\vert\leq n$
contribute; so Lemma 8.2; hence $\widehat{(\vert P_nu\vert^2)}(m)$ satisfies similar estimates
to $\widehat{(\vert u\vert^2)}(m)$, with the same constants.\par
\indent Let $\Phi :{\bf C}^{\ell -n}\rightarrow [0, \infty )$ be the convex function
$$\Phi (z_1, \dots ,z_{\ell -n})=\max_p\Bigl\{\Bigl\vert\sum_{t=n}^p z_{t-n}\Bigr\vert^4 : n\leq
p\leq \ell\Bigr\}\eqno(8.30)$$
\noindent associated with the fourth power of maximal partial sums. 
Then by the contraction principle from [12], the martingale maximal theorem in $L^4$ and
Khinchine's inequality we have
$$\eqalignno{\Bigl(\int_\Omega \Phi (d_n^{(m)}, \dots , d_\ell^{(m)})\mu
(du)\Bigr)^{1/4}&\leq 
\Bigl({\bf E}_{\varepsilon}\Phi (\delta_n\varepsilon_n, \dots , \delta_\ell
\varepsilon_\ell)\Bigr)^{1/4}\cr
&\leq {{4\sqrt{2}}\over{3}}\Bigl(\sum_{p=n}^\ell\delta_p^2\Bigr)^{1/2}\cr
&\leq {{4{\sqrt{2}}C_0K_2^2}\over{3\sqrt{\varepsilon}
n^{\varepsilon}}}.&(8.31)\cr}$$
\noindent The sequence $(\hat V(m))_{m\in {\bf Z}^2}$ is summable, so we deduce from (8.28) via the triangle
inequality in $L^2(\mu )$ and H\"older's inequality that
$$\eqalignno{&\Bigl(\int_\Omega \max_p\Bigl\{\vert U(P_pu)-U(P_nu)\vert^2 : n\leq p\leq 
\ell\Bigr\}\mu (du)\Bigr)^{1/2}&(8.32)\cr
&\quad \leq 
4\sum_{m\in {\bf Z}^2}\bigl\vert \hat V(m)\vert 
\Bigl( \int_\Omega \Phi (d_n^{(m)}, \dots , d_\ell^{(m)})\mu
(du)\Bigr)^{1/4}\cr
&\qquad\times\Bigl(\int_\Omega 
\bigl\vert \widehat{(\vert u\vert^2)}(m)\bigr\vert^4\mu (du)+
\int_\Omega \Phi (d_n^{(m)}, \dots , d_\ell^{(m)})\mu (du)\Bigr)^{1/4};\cr}$$
\noindent and hence by (8.31)  
$$\mu \Bigl\{ \max_p\bigl\{\vert U(P_pu)-U(P_nu)\vert^2 : n\leq p\leq \ell\bigr\} \geq \delta
\Bigr\}\rightarrow 0\qquad (\delta >0)\eqno(8.33)$$
\noindent as $\ell\geq n\rightarrow\infty$, so $U(P_nu)\rightarrow U(u)$ almost surely
and in $L^2(\mu )$ as $n\rightarrow\infty$. \par
\indent (ii)  We have
$${\hbox{Ent}}(\nu_n\mid \nu )=\int_{\Omega}\Bigl( U(P_nu)-U(u)+\log Z-\log
Z_n\Bigr)\nu_n(du),\eqno(8.34)$$
\noindent where the normalizing constants satisfy
$\lim\inf_{n\rightarrow\infty} Z_n\geq Z$, and the preceding arguments show that 
$\int_\Omega \vert U(P_nu)-U(u)\vert^2\mu (du)\rightarrow 0$ as
$n\rightarrow\infty$ and $\int_\Omega e^{2U(P_nu)}\mu (du)\leq C.$ 
Hence ${\hbox{Ent}}(\nu_n\mid \nu )\rightarrow 0$ as
$n\rightarrow\infty .$ By the transportation inequality (8.5), this implies ${\hbox{W}}_2(\nu_n, \nu
)\rightarrow 0$ as $n\rightarrow\infty$. Essentially, $\nu_n$ is the tensor
product of $\omega_n$ with a Gaussian measure on ${\hbox{H}}^{-s}$ with variance that 
converges to
zero as $n\rightarrow\infty$; indeed, the tail of the product (8.26) satisfies
$$\eqalignno{\int\sum_{m\in {\bf Z}^2; \vert m\vert\geq n}
{{a_m^2+b_m^2}\over{\vert m\vert^{2s}}}\prod_{m\in {\bf Z}^2; \vert m\vert\geq n}e^{-\vert
m\vert^2(a_m^2+b_m^2)/2} {{\vert m\vert^2da_mdb_m}\over{2\pi}}&=
\sum_{m\in {\bf Z}^2; \vert m\vert\geq n}{{2}\over{\vert m\vert^{2+2s}}}\cr
&\leq
{{4\pi}\over{s(n-1)^{2s}}}.\cr}$$
\noindent Hence ${\hbox{D}}_{L^2}((X^n,\Vert\cdot\Vert_{{\hbox{H}}^{-s}}, 
\omega_n),(\Omega_n,\Vert\cdot\Vert_{{\hbox{H}}^{-s}},\nu_n))\rightarrow 0$ as $n\rightarrow\infty$ as in
 [19, Example 3.8].\qquad{$\square$}\par
\vskip.05in
\indent Let
$\Delta={{\partial^2}/{\partial\theta^2_1}}+{{\partial^2}/
{\partial\theta_2^2}}$, and write
$$\Phi (u)(\theta , t)=\int_0^t  
 e^{i(t-\tau )\Delta}\Bigl( \bigl( \vert u\vert^2\ast
V\bigr)u\Bigr)(\theta , \tau
)\, d\tau .\eqno(8.35)$$
\noindent In Proposition 8.4, we verify that the solution of the G-P equation 
$$\eqalignno{ -i{{\partial u}\over{\partial t}}&=\Delta u+\bigl(\vert
u\vert^2\ast V\bigr) u,\cr
u(\theta , 0)&=\phi (\theta )&(8.36)\cr}$$
\noindent with $\phi\in\Omega \subset {\hbox{H}}^s({\bf T}^2; {\bf C})$ is given
by $u=u_0+w$, where $u_0(\theta , t)=e^{it\Delta}\phi (\theta )$ is the solution
of the free periodic Schr\"odinger equation with initial datum $\phi$ in the support of
Brownian loop on ${\hbox{H}}^{-s}$ and $w\in {\hbox{H}}^s$ is a
fixed point of $w\mapsto \Phi (u_0+w).$ \par 
\indent We say that $f:{\hbox{H}}^{-s}\rightarrow {\bf R}$ is a cylindrical function, if there
 exists a compactly supported
smooth function
$F:{\bf R}^n\rightarrow {\bf R}$ and $\xi_1, \dots ,\xi_n\in {\hbox{H}}^s$ such that
 $f(\phi )=
F(\langle \phi ,\xi_1\rangle , \dots \langle \phi ,\xi_n\rangle ).$ The following may be compared 
with Bourgain's results from [9, p. 132].\par
\par
\vskip.05in
\noindent {\bf Proposition 8.4} {\sl Let $0<s<1/68$, and let} $V\in
{\hbox{H}}^{\delta +2s+3/2}({\bf T}^2; {\bf R})$ {\sl for some $\delta
>0$ have $\hat V(0)=0$. Then for all $\eta>0$, there exists
$\Omega_\eta\subseteq \Omega$ and $L_\eta, t_\eta>0$ such that $\mu (\Omega_\eta )>1-\eta$ and \par
\indent (i) for all $\phi\in\Omega_\eta$ and 
 $u_0(\theta ,\tau )=e^{i\tau \Delta}\phi (\theta )$, 
the function} $\Phi (u_0)\in C([0,T]; {\hbox{H}}^s({\bf T}^2; {\bf C}))$ {\sl
for $T>0$ almost surely;\par
\indent (ii) $w\mapsto \Phi (u_0+w)$ is $L_\eta$-Lipschitz on
bounded subsets of} $C([0,T]; {\hbox{H}}^s({\bf T}^2; {\bf C}))$;\par
\indent {\sl (iii) the Cauchy problem (8.36) has a solution
$u(\theta ,  t)$ for 
$t\in [0, t_\eta ]$ for all $\phi \in\Omega_\eta$;\par
\indent (iv) $\phi (\theta )  \mapsto u(\theta, t)$ for
$\phi\in\Omega_\eta$ induces a measure on} ${\hbox{H}}^{-s}$ {\sl which 
satisfies the $T_1$ transportation inequality, and is invariant in the sense
 that all cylindrical functions
satisfy}
$$\int_{\Omega_\eta} f(u(\,\cdot \, ,t))\nu (d\phi )= 
\int_{\Omega_\eta} f(\phi )\nu (d\phi )\qquad (0\leq t<t_\eta ).\eqno(8.37)$$
\vskip.05in
\noindent {\bf Proof.} (i) We write $\Vert a\Vert_*=1+\Vert a\Vert$. Note that
$(\Omega, \mu )$ is invariant under the operation $\phi (\theta )\mapsto
e^{i\tau \Delta}\phi (\theta ).$  The
integral (8.35) may be expressed in Fourier coefficients as
$$\eqalignno{\Phi &(u_0)(\theta ,t)&(8.38)\cr
&=\sum_{m\in {\bf Z}^2\setminus \{ 0\}} \Bigl[ \sum_{j,k: j+k=m}\hat \phi
(j)\overline{ \hat \phi (-k)} {{e^{it(\Vert \ell\Vert^2-\Vert\ell
+m\Vert^2)} -1}\over{i( \Vert \ell\Vert^2-\Vert\ell+m\Vert^2+\Vert
j\Vert^2-\Vert k\Vert^2)}} \Bigr] \hat V(m)\sum_\ell e^{i(\ell
+m)\cdot\theta}\hat \phi (\ell ),\cr}$$
\noindent and we split this sum into four
cases, according to the values of $j$ and $k$ in the inner sum, and
then according to $\ell$ and $m$ in the outer sums. First we note that
in the inner sum in square brackets $\Vert j\Vert^2-\Vert
k\Vert^2=(2j-m)\cdot m$, so 
we split the index set as $\{ (j,k)\in {\bf Z}^2\times {\bf Z}^2:
j+k=m\} =G(\ell ,m)\sqcup B(\ell ,m\}$ where 
$$G(\ell ,m)=\Bigl\{ (j,k):
j+k=m; \bigl\vert \Vert \ell\Vert^2-\Vert
m+\ell\Vert^2+(2j-m).m\bigr\vert\geq 2^{-2} \bigl\vert \Vert
\ell\Vert^2-\Vert \ell +m\Vert^2\bigr\vert \Bigr\},\eqno(8.39)$$
\noindent and the complementary set
$$B(\ell ,m)=\Bigl\{ (j,k):
j+k=m; \bigl\vert \Vert \ell\Vert^2-\Vert
m+\ell\Vert^2+(2j-m).m\bigr\vert<2^{-2} \bigl\vert \Vert
\ell\Vert^2-\Vert \ell +m\Vert^2\bigr\vert \Bigr\},\eqno(8.40)$$
\noindent so that $B(\ell ,m)$ is the set of integral lattice points in a strip in ${\bf
R}^2$ which has axis perpendicular to $m$ and width $ \bigl\vert \Vert
\ell\Vert^2-\Vert \ell +m\Vert^2\bigr\vert$. Now the sum
$$\sum_{(j,k)\in G(\ell ,m)} {{\hat \phi (j)\overline{\hat \phi
(-k)}(\Vert \ell \Vert^2-\Vert \ell +m\Vert^2)}\over{\Vert
\ell\Vert^2-\Vert\ell+m\Vert^2+\Vert j\Vert^2-\Vert k\Vert^2}}\eqno(8.41)$$
\noindent is exponentially square integrable by Lemma 8.2. Then we take
the complementary contribution to the inner sum of (8.38) to be 
$$\eqalignno{ \Bigl\vert &\sum_{(j,k)\in B(\ell ,m)} 
{{\hat \phi (j)\overline{\hat \phi
(-k)}}\over{1+\bigl\vert \Vert
\ell\Vert^2-\Vert\ell+m\Vert^2+\Vert j\Vert^2-\Vert
k\Vert^2\bigr\vert}}\Bigr\vert\cr
&\leq \sum_{(j,k)\in B(\ell ,m)}{{K_2^2}\over{ \Vert
j\Vert^{\varepsilon+3/4}_*\Vert k\Vert^{\varepsilon +3/4}_*\bigl\vert
\Vert \ell\Vert^2-\Vert\ell+m\Vert^2+\Vert j\Vert^2-\Vert
k\Vert^2\bigr\vert}}\cr
&\leq {{K_2^2}\over{\bigl\vert \Vert \ell\Vert^2-\Vert \ell
+m\Vert^2\bigr\vert^{1/16}}}\sum_{(j,k)\in B(\ell ,m)}
\Bigl\{ {{1}\over{\Vert j\Vert^{\varepsilon+1/2}_*\Vert 
m-j\Vert^{\varepsilon +1/2}}}\Bigr\}\cr
&\qquad\times \Bigl({{1}\over{ \Vert
j\Vert^{1/8}_*\Vert m-j\Vert^{1/8}_* \bigl\vert \Vert \ell\Vert^2-\Vert
\ell +m\Vert^2+2j\cdot m-\Vert m\Vert^2\bigr\vert^{15/16}}}\Bigr).&(8.42)\cr}$$
Then we split $j=j_\perp +j_m$, where $j_\perp$ is perpendicular to
$m$, and $j_m$ parallel to $m$; the sum in braces is dominated by the
corresponding sum over $j_\perp$ and is bounded; while the sum in
parentheses is dominated by the corresponding sum over $j_m$ and is
also bounded; so the whole expression (8.42) is 
$$\leq C{{K_2^2}\over{\bigl\vert \Vert \ell\Vert^2-\Vert \ell
+m\Vert^2\bigr\vert^{1/16}}}.\eqno(8.43)$$
\indent We deduce that for all $\eta>0$, there exist a subset
$\Omega_\eta\subset \Omega$ with $\mu (\Omega_\eta )>1-\eta$ and a constant
$C_\eta$ such that 
$$\eqalignno{\bigl\Vert \Phi (u_0)\bigr\Vert_{{\hbox{H}}^s}&\leq
C_\eta\sum_m\vert \hat V(m)\vert \Bigl\Vert \sum_\ell {{e^{i(m+\ell
)\cdot\theta} \hat \phi (\ell )}\over{\bigl\vert \Vert\ell\Vert^2-\Vert
\ell +m\Vert^2\bigr\vert^{1/16}}}\Bigr\Vert_{{\hbox{H}}^s}\cr
&\leq C_\eta\sum_m\Vert \hat V(m)\vert\Vert m\Vert^{2s}_* \Bigl[\sum_\ell
\Bigl( {{\Vert \ell +m\Vert^{2s}\Vert\ell\Vert_*^{2s}}\over{\Vert
m\Vert^{2s}_*\bigl\vert \Vert
\ell\Vert^2-\Vert\ell+m\Vert^2\bigr\vert^{1/8}}}\Bigr)
{{\vert\hat\phi(\ell
)\vert^2}\over{\Vert\ell\Vert_*^{2s}}}\Bigr]^{1/2}.&(8.44)\cr}$$
\noindent We split this sum into a sum over the index set 
$$A=\Bigl\{ (\ell ,m )\in {\bf Z}^2\times {\bf Z}^2: \bigl\vert \Vert
\ell\Vert^2-\Vert \ell+m\Vert^2\big\vert^{1/8}\geq
\Vert\ell\Vert_*^{4s}\Bigr\}$$
\noindent and a sum over the
complementary set $A^c$. On $A$, the factor in parentheses from (8.44) is bounded, so the upper
bound $\sum_m\vert \hat V(m)\vert \Vert m\Vert^{2s}\Vert \phi\Vert_{{\hbox{H}}^{-s}}$
is immediate. On $A^c$, we use the bound $\vert\hat \phi (j)\vert\leq K_2\Vert j\Vert_*^{-\varepsilon
-3/4}$, and for each $m$, we compare the sum over $(\ell
,m)\in A^c$ with an integral in polar coordinates $(r, \psi )$ over the  
region
$$\bigl\{ (r, \psi )\in (1, \infty )\times (-\pi ,\pi ): 2\Vert
m\Vert r\vert\sin\psi\vert \leq \Vert m\Vert^2+r^{32s}\bigr\};\eqno(8.45)$$
so we have a bound on $\sum_{\ell\in A^c}$ of 
$$\eqalignno{\sum_{\ell :\vert \Vert \ell\Vert^2-\Vert \ell +m\Vert^2\vert
<\Vert\ell\Vert^{32s}} &{{K_2^2}\over{ \Vert m\Vert_*^{2s}\vert \Vert
\ell \Vert^2-\Vert m+\ell \Vert^2\vert^{1/8}\Vert
\ell\Vert_*^{2\varepsilon+3/2}}}\cr
&\leq 2K_2^2\int_1^\infty
r^{2s-3/2-2\varepsilon }\int_0^{(r^{32s}+\Vert m\Vert^2)/2r\Vert m\Vert} d\psi \,
rdr\cr
&\leq 2CK_2\Bigl(\Vert m\Vert+{{1}\over{\Vert m\Vert}}\Bigr).&(8.46)\cr}$$
\noindent The series $\sum_m \vert \hat V(m)\vert \Vert m\Vert^{2s+1/2}$ converges, 
so $\Phi (u_0)$ belongs to $C([0,T]; {\hbox{H}}^s).$\par
\indent (ii) In this proof, we use concentration of measure to prove
Lipschitz continuity of a function; this reverses the usual flow of the
theory as in [5, 21]. For $v$ and $w$ in the unit ball of\par
\noindent $C([0, T]; {\hbox{H}}^s({\bf
T}^2; {\bf C}))$, We have 
$$\eqalignno{\Phi (v+u_0)&(\theta ,t)-\Phi (w+u_0)(\theta ,t)&(8.47)\cr
&=\int_0^t e^{i(t-\tau )\Delta}
\Bigl( \bigl( \vert u_0\vert^2\ast V\bigr) (v-w)\Bigr) (\theta ,\tau )\,
d\tau \cr
&\quad +\int_0^t e^{i(t-\tau )\Delta}
\Bigl( \bigl( (\vert v\vert^2+\bar vu_0+v\bar u_0)\ast V\bigr) (v-w)\Bigr) (\theta ,\tau )\,
d\tau \cr
&\quad+\int_0^t e^{i(t-\tau )\Delta}
\Bigl( \bigl( ( (v-w)\bar v +w(\bar v-\bar w)+u_0(\bar v-\bar w)+\bar
u_0(v-w))\ast V\bigr) w\Bigr) (\theta ,\tau )\,
d\tau . \cr}$$
\noindent In the final integral, we can use the simple bound 
$$\bigl\vert \widehat{u_0(\bar v-\bar w)}(m)\bigr\vert\leq \bigl\Vert
u_0\bigr\Vert_{{\hbox{H}}^{-s}}\bigl\Vert
v-w\bigr\Vert_{{\hbox{H}}^{s}}\leq K_1 \bigl\Vert
v-w\bigr\Vert_{{\hbox{H}}^{s}},\eqno(8.48)$$
\noindent and similar bounds on the other terms; the terms in the middle
integral are treated similarly. The first integral, we use the
probabilistic estimate of Lemma 8.2: for all $\eta>0$ there exist $L_\eta>0$
and a subset $\Omega_\eta\subseteq \Omega$ such that $\mu (\Omega_\eta
)>1-\eta$ and 
$$\sum_m\bigl\vert \widehat{(\vert u_0\vert^2)}(m)\bigr\vert \vert \hat
V(m)\vert \Vert m\Vert^{2s}\leq L_\eta\qquad (u_0(\theta , 0)\in \Omega_\eta
) ,$$
\noindent so there exists $C>0$ such that
$$\sup_{0<t<T}\Vert \Phi (u_0+v)(\theta ,t)-\Phi
(u_0+w)(\theta ,t)\Vert_{{\hbox{H}}^{s}}\leq
CT(1+L_\eta)\sup_{0<\tau<T}\Vert v(\theta ,\tau)-w(\theta
,\tau)\Vert_{{\hbox{H}}^{s}}.\eqno(8.49)$$  
\indent (iii) By (i), we have $T>0$ such that $K_0=\sup_{0<t<T}\Vert
\Phi(u_0)(\theta ,t)\Vert_{{\hbox{H}}^s}$ is finite for all $\phi\in
\Omega_\eta.$ Now by (8.49) we can shrink the time interval to $[0, t_\eta]$ where
$0<t_\eta<T$, and ensure that 
$$B_\eta =\Bigl\{ w\in C([0, t_\eta ];
{\hbox{H}}^s({\bf T}^2; {\bf C})); \sup_{0<t<t_\eta}\Vert
w(\theta ,t)\Vert_{{\hbox{H}}^s}\leq 2K_0\Bigr\}\eqno(8.50)$$
\noindent contains $\Phi (u_0)$ and $w\mapsto \Phi (u_0+w)$ is $(1/2)$-Lipschitz
on $B_\eta$. Indeed, we have
$$\eqalignno{\sup_{0<t<t_\eta}\Vert
\Phi(u_0+w)(\theta ,t)\Vert_{{\hbox{H}}^s}&\leq \sup_{0<t<t_\eta}\Vert
\Phi(u_0+w)(\theta ,t)-\Phi (u_0)(\theta
,t)\Vert_{{\hbox{H}}^s}+\sup_{0<t<t_\eta}\Vert
\Phi(u_0)(\theta ,t)\Vert_{{\hbox{H}}^s}\cr
&\leq 2^{-1}\sup_{0<t<t_\eta}\Vert
w(\theta ,t)\Vert_{{\hbox{H}}^s}+K_0\cr
&\leq 2K_0.&(8.51)\cr}$$
\noindent By Banach's fixed point theorem, there exists $w\in B_\eta$ such that
$w=\Phi (u_0+w);$ thus we obtain a solution $u(\theta ,t)=u_0(\theta
,t)+w(\theta ,t)$ of G-P (8.36) for
$0<t<t_\eta$.\par
\indent (iv) We do not assert that $\phi \mapsto \Phi (u_0+v)$ is Lipschitz; hence we need an indirect proof of (iv) instead of deducing it from Theorem
8.1. The fixed point $w$ satisfies $\Vert
w(\,\cdot\,  ,t)\Vert_{{\hbox{H}}^s}\leq 
2\Vert \Phi (u_0)(\, \cdot\, ,t)\Vert_{{\hbox{H}}^s},$ hence 
$$\Vert u(\, \cdot\, ,t)\Vert_{{\hbox{H}}^{-s}}\leq
\Vert\phi\Vert_{{\hbox{H}}^{-s}}
+\Vert \Phi (u_0)(\, \cdot\, ,t)\Vert_{{\hbox{H}}^s}\eqno(8.52)$$
\noindent so there exists $\kappa>0$ such that
$$\int_{\Omega_\eta} \exp \bigl( \kappa \Vert u(\, \cdot\, ,t)
\Vert^2_{{\hbox{H}}^{-s}}\bigr) \nu (d\phi )\eqno(8.53)$$
\noindent is finite. Hence the measure induced on ${\hbox{H}}^{-s}$ from $\mu$
on $\Omega_\eta$  by
$\phi\mapsto u(\, \cdot\, ,t)$ satisfies a $T_1$ transportation inequality by Bobkov and G\"otze's
criterion, as in [21, Theorem 22.10]. \par 
\indent  Let $u_n$ be the solution
of the $GP$ equation with finite-dimensional Hamiltonian $H_n$ as in (7.3) and initial
data $\phi_n(\theta) =\sum_{k: 0<\vert k\vert \leq n} e^{ik\cdot\theta
}(\gamma_k+i\tilde\gamma_k)/\Vert k\Vert$, and we regard $u_n(\, \cdot \, ,t)$ as a random
variable for $\phi\in\Omega_\eta$. We have
$$\eqalignno{\Vert & u(\, \cdot \, ,t)-u_n(\, \cdot \, ,t)\Vert_{{\hbox{H}}^{-s}}\cr
&\leq 2\Vert
\phi-\phi_n\Vert_{{\hbox{H}}^{-s}}+2\Bigl\Vert\int_0^te^{i(t-\tau )\Delta }\Bigl(\bigl(\vert w_n+e^{i\tau\Delta}\phi_n\vert^2\ast
V\bigr)\bigl(e^{i\tau\Delta }\phi_n-e^{i\tau\Delta}\phi\bigr) \Bigr) d\tau 
\Bigr\Vert_{{\hbox{H}}^{-s}}\cr
&\quad +2\Bigl\Vert\int_0^te^{i(t-\tau )\Delta }\Bigl( \bigl(
\vert w_n+e^{i\tau\Delta}\phi\vert^2-\vert w_n+e^{i\tau\Delta}\phi_n\vert^2\bigr)\ast
V\bigr)\bigl(w+e^{i\tau\Delta }\phi \bigr)\Bigr) d\tau \Bigr\Vert_{{\hbox{H}}^{-s}}
.&(8.54)\cr}$$

As in (iii), one can show that $u_n$ converges to
$u$ in the sense that
$$\int_{\Omega_\eta} \Vert u_n(\,\cdot \, ,t)-u(\, \cdot \,
,t)\Vert^2_{{\hbox{H}}^{-s}}\mu (d\phi)\rightarrow 0\eqno(8.55)$$
\noindent as $n\rightarrow\infty. $ By Liouville's theorem applied to $H_n$, the corresponding Gibbs measure on
phase space is invariant under the flow generated by the canonical equations of motion. Hence by
Corollary 8.3, we have weak convergence of the Gibbs measures, so
$$\eqalignno{ \int_{\Omega_\eta} f(u(\,\cdot \, ,t))\nu (d\phi )&=\lim_{n\rightarrow\infty}
\int_{\Omega_\eta} f(u_n(\,\cdot \, ,t))\nu_n (d\phi )\cr
&=\lim_{n\rightarrow\infty}\int_{\Omega_\eta} f(\phi )\nu_n (d\phi )\cr
&=\int_{\Omega_\eta} f(\phi )\nu (d\phi ).&(8.56)\cr}$$
\rightline{$\square$}\par
\vskip.05in    
\noindent {\bf Acknowledgement.} The author thanks Graham Jameson and Daniel
Elton for providing some
estimates.\par 
\vskip.05in
\noindent {\bf References}\par
\noindent [1] S. Aida and I. Shigekawa, Logarithmic Sobolev inequalities
and spectral gaps: perturbation theory, 
{\sl J. Funct. Anal.} {\bf 126} (1994), 448--475.\par
\noindent [2] S.G. Bobkov and M. Ledoux, From Brunn--Minkowski to Brascamp--Lieb
and to logarithmic Sobolev inequalities, {\sl Geom. Funct. Anal.} {\bf 10} (2000),
1028--1052.\par 
\noindent [3] G. Blower, A logarithmic Sobolev inequality for the invariant measure of the periodic
Korteweg--de Vries equation, {\sl Stochastics} {\bf 84} (2012), 533--542.\par
\noindent [4] G. Blower and F. Bolley, Concentration of measure on
product spaces with applications to Markov processes, {\sl Studia
Math.} {\bf 175} (2006), 47--72.\par
\noindent [5] S.G. Bobkov, I. Gentil and M. Ledoux, Hypercontractivity of
Hamilton--Jacobi equations, {\sl J. Math.
Pures Appl. (9)} {\bf 80} (2001), 669--696.\par
\noindent [6] J. Bourgain, Periodic nonlinear Schr\"odinger equation and
invariant measures, {\sl Comm. Math. Phys.} {\bf 166} (1994), 1--26.\par
\noindent [7] J. Bourgain, On the Cauchy and invariant measure problem for the periodic
Zakharov system, {\sl Duke Math. J.} {\bf 76} (1994), 175--202.\par
\noindent [8] J. Bourgain, Invariant measures for the Gross--Piatevskii
equation, {\sl J. Math. Pures Appl. (9)} {\bf 76} (1997), 649--702.\par
\noindent [9] J. Bourgain, {\sl Global Solutions of Nonlinear Schr\"odinger
Equations}, (American Mathematical Society, Providence RI, 1999).\par 
\noindent [10] D.C. Brydges and G. Slade, Statistical mechanics of the
 2-dimensional focusing
nonlinear Schr\"odinger equation, {\sl Comm. Math. Phys.} {\bf 182} (1996), 485--504.\par
\noindent [11] R. Holley and D. Stroock, Logarithmic Sobolev inequalities and stochastic
Ising models, {\sl J. Statist. Physics} {\bf 46} (1987), 1159-1194.\par
\noindent [12] J. Jakubowski and S. Kwapien, On multiplicative systems of
functions, {\sl Bull. Acad. Polon. Sci. Ser. Sci. Math.} {\bf 27} (1979),
689--694.\par  
\noindent [13] K. Kirkpatrick, Solitons and Gibbs measures for nonlinear Schr\"odinger 
equations, {\sl Math. Model. Nat. Phenom.} {\bf 7} (2012), 95--112.\par
\noindent [14] J.L. Lebowitz, Ph. Mounaix and W.-M. Wang, Approach to equilibrium for stochastic NLS,
{\sl Commun. Math. Phys.} {\bf 321} (2013), 69-84.\par
\noindent [15] J.L. Lebowitz, H.A. Rose and E.R. Speer, Statistical 
mechanics of the nonlinear\par
\noindent Schr\"odinger equation, {\sl J. Statist. Phys.} {\bf 50} (1988), 657-687.\par
\noindent [16] H. P. McKean, Statistical mechanics of nonlinear wave equations 
IV: cubic Schr\"odinger, {\sl Comm. Math. Phys.} {\bf 168} (1995), 479--491.\par 
\noindent [17] Y. Pomeau and S. Rica, Dynamics of a model of a supersolid, {\sl Phys. Rev. Lett.} 
{\bf 72} (1994), 2426--2430.\par
\noindent [18] L. Schwartz, {\sl Radon measures on arbitrary topologival
spaces and cylindrical measures}, (Oxford University Press, 1973).\par 
\noindent [19] K.-Th. Sturm, On the geometry of metric measure spaces I,
{\sl Acta Math.} {\bf 196} (2006), 65--131.\par
\noindent [20] T. Tao and V. Vu, Random matrices: universality of local
spectral statistics of non-Hermitian matrices, {\sl Ann. Probab.} {\bf 43}
(2015), 782--874.\par
\noindent [21] C. Villani, {\sl Optimal transport: old and new},
(Springer-Verlag, Berlin, 2009).\par
\vfill
\eject
\end   
We introduce the norm $\Vert\, \cdot\Vert$ on $X^n$ from the embedding in
$L^p$ associated with the trigonometric system, so
$$\Vert (x_k)_k\Vert =\Bigl(\int_{{\bf T}^D}\bigl\vert \sum_{k} x_ke^{ik\cdot\theta}\bigr\vert^p
{{d^D\theta}\over{(2\pi)^D}}\Bigr)^{1/p}\eqno(2.11)$$
\noindent with the dual norm $\Vert\, \cdot\Vert_*$ for the sesquilinear pairing
$\sum_{k}x_k\bar y_k$. By M. Riesz's theorem, for $1<p<\infty$ there exist constants 
depending on $p$ and $D$ but not on $n$ such that
$$\Vert (y_k)_k\Vert_* \sim \Bigl(\int_{{\bf T}^D}\bigl\vert \sum_{k}
y_ke^{ik\cdot\theta}\bigr\vert^{p^*}
{{d^D\theta}\over{(2\pi)^D}}\Bigr)^{1/p^*}\eqno(2.12)$$
\noindent where $p^*=p/(p-1)$. The norm of $L^p$ is uniformly convex and uniformly
smooth, and the unique geodesic connecting two points is the line segment between
them. Also, the norm $\Vert\cdot \Vert$ is equivalent to the Euclidean norm on $X^n$, with constants
depending upon $n$.\par

\noindent {\bf 7. Concentration and transportation inequalities for $D=1$} 
\vskip.05in 
\noindent In this section, we focus on $D=1$ and deduce some concentration and
transportation inequalities from previous sections. For $1<\beta\leq 2$ and $0<a<\infty$, we introduce 
$$L(x)=\cases{ \beta^{-1}x^\beta,& for $0\leq x<a$;\cr 
      a^{\beta -1}(x-a)+\beta^{-1}a^\beta, & for $a\leq x$,\cr}\eqno(7.1)$$
\noindent so that $L:[0, \infty )\rightarrow [0, \infty )$ is strictly increasing and convex
with $L(0)=0$, and has Legendre transform      
$$L^*(x)=\cases{ (1-\beta^{-1})x^{\beta /(\beta -1)},& for $0\leq x\leq a^{\beta -1}$;\cr 
      \infty, & for $a^{\beta -1}<x$.\cr}\eqno(7.2)$$
\noindent For a bounded and continuous $F:X^n\rightarrow {\bf R}$, we define for
$t>0$
$$Q_tF(x)=\inf\Bigl\{ F(y)+tL\bigl( t^{-1}\Vert x-y\Vert\bigr): y\in X^n\Bigr\}\qquad (x\in
X^n).\eqno(7.3)$$
\noindent The map $F\mapsto Q_tF$ is monotone and concave, although nonlinear, for the
pointwise operations on $F$. We have chosen $L$ with arbitrary $0<a<\infty$ so that all the $Q_tF(x)$ for all
$t>0$ are Lipschitz. By slightly adapting results of [22, Theorem 22.46] to the geodesic space $(X^n,
\Vert\cdot \Vert )$, one can check that $Q_t(Q_sF))=Q_{t+s}F$
for Lipschitz $F$, and 
$${{\partial}\over{\partial t}}Q_tF(x)=-L^*\bigl( \Vert \nabla
Q_tF\Vert_*\bigr).\eqno(7.4)$$
\noindent The logarithmic moment generating function of $\mu$  is $\Lambda (f)=\log\int
e^{\langle f, \phi \rangle} \nu (d\phi )$. Taking $F(\phi )=\langle f,\phi \rangle$, we can
express the logarithmic moment generating function of $\nu$ in terms of $Q_tF$
by
$$\Lambda (\alpha f)=\log\int e^{\alpha Q_tF(\phi
)}\nu (d\phi )+t\alpha L^*(\Vert f\Vert_*)\qquad (\alpha, t>0).\eqno(7.5)$$
\noindent More generally, for $F$ convex, $Q_tF(x)\geq F(x)-tL^*( \Vert\nabla F(x)\Vert_*).$\par
\vskip.05in
\noindent {\bf Definition.} Say that $\nu\in {\hbox{Prob}}_2(X^n)$ satisfies the
defective generalized logarithmic Sobolev
inequality, if there exist $\alpha>0$ and $\gamma\geq 0$ such that 
$${\hbox{Ent}}(\omega \mid \nu )\leq {{1}\over{\alpha}}\int_XL^*(\Vert \nabla
U\Vert_*)\,d\omega+\gamma\eqno(7.6) $$ 
\noindent holds for all $\omega\in {\hbox{Prob}}_2(X^n)$ of the form $\omega =e^U\nu$ where 
$U:X\rightarrow {\bf R}$ is Lipschitz and satisfies 
$\Vert \nabla U\Vert_*\leq a^{\beta -1}$. If $\gamma =0$, then $\nu$ satisfies a
generalized logarithmic Sobolev inequality. If (7.6) holds with $\gamma =0$, $\beta =2$ and some
$\alpha>0$ for all 
$a\rightarrow\infty$, we say that $\nu$ satisfies $LSI(\alpha )$.
Simple example show that the
defective logarithmic Sobolev inequality with $\gamma >0$ does not imply 
$LSI (\alpha )$ in
general. See [1] and [22, section 22] for discussion.\par
\vskip.05in

\noindent {\bf Proposition 7.1} {\sl Let $D=1$ and $p=4$, and let $\nu_N$ be the Gibbs measure on $\Omega_N$
associated with $H_4$ as in (1.1).\par
\indent (i)  For all $N$ 
there exist $\lambda, \alpha>0$ such that $\nu_N$ satisfies $LSI(\alpha )$.\par 
\indent (ii) The concentration inequality  
$$\int_{\Omega_N}\exp\bigl( \alpha Q_1F(x)\bigr)\nu_N (dx)\leq \exp\Bigl(
 \alpha\int_{\Omega_N} 
F(x)\,\nu_N (dx)\Bigr)\eqno(7.7)$$
\noindent and the Poincar\'e (spectral gap) inequality
$$\int_{\Omega_N}\Bigl( F(u)-\int Fd\nu_N\Bigr)^2\, \nu_N(du)\leq
{{2}\over{\alpha}}\int_{\Omega_N}L^*\bigl(\Vert\nabla F(u)\Vert_{L^2}\bigr)\nu_N
(du)\eqno(7.8)$$
\noindent with $\beta =2$ hold for all Lipschitz $F:(\Omega_N, L^2)\rightarrow {\bf R}$.}\par
\vskip.05in
\noindent {\bf Proof.} (i) See Proposition 3.1\par
\indent (iv)  Suppose that (7.6) holds for some $\beta>1$ and $\gamma =0$. Then 
$$t^{1-\beta}\log\Bigl( \int e^{\alpha t^{\beta
-1}Q_tF}d\nu\Bigr)\eqno(7.9)$$
\noindent is monotone decreasing on $(0, 1)$ and converges to $\alpha \int Fd\nu$ as
$t\rightarrow 0+$. For $\beta =\beta/ (\beta-1)=2$, we have $L^*(st)\leq t^2L^*(s)$ for
all $0<t<1$ and $0<s$. The concentration inequality and spectral gap inequality follow from 
the generalized logarithmic Sobolev
inequality for all Lipschitz $F$ by Theorem 22.28 of [22]. Given $F$, 
we can choose $a>\Vert
F\Vert_{Lip}$ in the definition of $L$, and thereby ensure that all the quantities in the
integrals are finite, and $a$ does not appear in the statement. See also the
exponential integrability results of Hebisch and Zegarlinski [13].\par
\rightline{$\square$}\par
\vskip.05in

\noindent {\bf Proposition 7.2} {\sl Let $D=1$ and $p=6$, and let $\nu_N$ be the Gibbs measure on $\Omega_N$
associated with $H_6$ as in (1.1). There exists $N_0>0$ such that for all $N<N_0$ where exist 
$a,\lambda, \alpha >0$ such that $\nu_N$ satisfies a defective generalized logarithmic Sobolev
inequality with $\beta =6$.}\par

\vskip.05in
\noindent {\bf Proof.} The proof is suggested by perturbation results of 
Aida and Shigekawa [1]. The potential
$$U_p(u)={{1}\over{p}}\int_{\bf T} \vert u(\theta )\vert^p \, d\theta\eqno(7.10)$$
\noindent has gradient, for $u,w\in X^n$, 
$$\eqalignno{ \langle \nabla U_p(u), w\rangle_{L^2}&=\lim_{t\rightarrow 0}
t^{-1} \bigl( U_p(u+tw)-U_p(u)\bigr)\cr
&={{1}\over{2}}\int_{\bf T}\vert u(\theta )\vert^{p-2}\bigl( \bar u(\theta
)w(\theta )+u(\theta )\bar w(\theta )\bigr)\, d\theta &(7.11)\cr}$$
\noindent so by the Cauchy--Schwarz inequality,
$$\bigl\Vert\nabla U_p(u)\bigr\Vert^2_{L^2}\leq \int_{\bf T}\vert u(\theta
)\vert^{2(p-1)}\, d\theta.\eqno(7.12)$$
\indent Let $\mu$ be the probability measure on $\Omega_N$ that is induced from
Wiener loop by random Fourier series by the map $\omega\mapsto b_\omega (\theta)\in L^2({\bf T})$
as in (1.3). Then by results of Federbush and Gross,
$\mu$ satisfies $LSI(\alpha_0)$ for some $\alpha_0>0$. The Gibbs measure 
for the Hamiltonian may be expressed
as
$$\nu_N(du)={\bf I}_{\Omega_N}(u)e^{U_p(u)}d\mu (u),\eqno(7.13)$$
 then let
$a=(5/(3\alpha_0 ))^{1/4}$ so that $L^*(x)\geq x^2/2$ for all $x\geq 0$. Hence $\mu$ 
satisfies
$$\int g\log g\, d\mu \leq {{2}\over{\alpha_0}}\int L^*(\Vert \nabla \log g\Vert
)gd\mu\eqno(7.14)$$
\noindent so we introduce the probability measure $d\nu =e^{U}d\mu$ and substitute $g=e^{U}f$, where $f$ is a probability
density function with respect to the $\nu$ with $\log f$ Lipschitz. Hence we have, after a
little rearranging
$$\int f\log f d\nu \leq {{2^{7/6}}\over{\alpha_0}}\int L^*(\Vert \nabla \log f\Vert ) fd\nu
+{{2^{7/6}}\over{\alpha_0}}\int L^* (\Vert\nabla U\Vert )fd\nu +\int
(-U)fd\nu .\eqno(7.15)$$
\noindent Now we use Young's inequality $hf\leq (1/2)f\log f +(1/2)e^{2h-1}$ on the final terms on the
right-hand side, obtaining  
$${{1}\over{2}}\int f\log f d\nu \leq 
{{2^{7/6}}\over{\alpha_0}}\int L^*(\Vert \nabla \log f\Vert )
fd\nu+{{1}\over{2}}\int\exp\Bigl( {{2^{13/6}}\over{\alpha_0}} L^*(\Vert \nabla U\Vert
)-2U-1\Bigr) d\nu.\eqno(7.16)$$
\noindent This has the form of a defective generalized logarithmic Sobolev inequality,
provided that the final integral is finite. Now for $N<N_0$ and $U(u)=\lambda\int u^6 dx$,
we have $L^*(\Vert \nabla U\Vert )=\lambda^{6/5}\int u^6$, so when $\lambda >0$
is sufficiently small, the main theorem of [18] shows that (7.16) is finite.\par 
\rightline{$\square$}\par
\vskip.05in
\indent Proposition 7.1 (ii) leads to the
following interpretation of relative entropy as a rate function, in the style of
Sanov. For $2\leq p<6$,
consider the metric probability space $(\Omega_N, \Vert\, \cdot\,\Vert_{L^p},
\nu_N)$. Let ${\bf X}={\hbox{Prob}}_p(\Omega )$ have the
Wasserstein metric ${\hbox{W}}_p$ for cost function $\Vert f-g\Vert^p_{L^p}$, so $({\bf X},
{\hbox{W}}_p)$ is a complete and separable metric space. Given mutually
independent $f_j\in \Omega_N$, each distributed according to $\nu_N$, we let
$f=(f_1, \dots ,f_n)$ and introduce the
empirical measure 
$$\varepsilon^f_n={{1}\over{n}}\sum_{j=1}^n \delta_{f_j}\eqno(7.17)$$
\noindent so $\varepsilon^f_n$ is a random element of ${\bf X}$ with distribution
determined by the product measure $\nu_N^{\otimes n}$. For an open subset
$B$ of $\Omega_N$,  $\varepsilon_n(B)$ is the relative frequency with
which $f_1, \dots ,f_n$ is in $B$.\par
\vskip.05in
\noindent {\bf Proposition 7.3} {\sl For $2\leq p<6$ let $\nu_N$ be the
Gibbs measure and let}
 $E(\omega )={\hbox{Ent}}(\omega \mid \nu_N)$.\par
\indent {\sl (i) Then $E$ is a good rate function; so that 
$E$ is lower semi-continuous, and for all $t>0$, the level subset $\{ \omega : E(\omega )\leq t\}$ of}      
$({\bf X}, {\hbox{W}}_p)$ {\sl is compact and convex.\par
\indent (ii) There exists $\alpha >0$ such that} $\{ \omega :
 E(\omega )\leq t\}\subseteq \{ \omega : {\hbox{W}}_p(\omega ,
\nu_N)\leq \sqrt{2t/\alpha}\}$ {\sl for all $t>0$.\par
\indent (iii) For all open subsets $G$ of} $(E, {\hbox{W}}_p)$
$$-\inf_\omega \{ E(\omega ): \omega\in G\}\leq \lim\inf_{n\rightarrow\infty}
 {{1}\over{n}}\log\nu_N^{\otimes n}(\{ f\in \Omega_N^n: \varepsilon^f_n\in
G\}).\eqno(7.18)$$
\indent {\sl (iv) For all closed subsets $K$ of} $(E, {\hbox{W}}_p)$
$$\lim\sup_{n\rightarrow\infty} {{1}\over{n}}\log\nu_N^{\otimes
n}(\{ f\in \Omega_N^n:\varepsilon^f_n\in K\})\leq
-\inf_\omega \{
E(\omega ): \omega\in K\}.\eqno(7.19)$$
\vskip.05in
\noindent {\bf Proof.} (i) This is standard [12, p. 36].\par
\indent (ii)  This is a reformulation of the transportation inequality (2.3) 
defining $T_p$.\par
\indent (iii) and (iv) The map $f\mapsto \varepsilon_n^f$ is
$1/\sqrt{n}$-Lipschitz $(\Omega_N^n, (L^2)^n)\rightarrow
({\hbox{Prob}}_2(\Omega_N), {\hbox{W}}_2)$, and hence is continuous. Consequently,  
$G\mapsto \nu_N^{\otimes n}(\{ f\in \Omega_N^n: \varepsilon^f_n\in
G\})$ defines a measure on ${\bf X}$. Observe that $\int_{\Omega_N} \exp( {\rho \Vert
f\Vert_{L^p}^p})\nu_N(df)$ is finite for all $\rho>0$. Hence we fall within the
scope of Theorem 1.1 of [23].\par
\rightline{$\square$}\par
\indent The cases $s=1$ and $s=2$ are of particular importance, as noted by Marton and Talagrand. 
The inequality $LSI(\alpha )$ implies $T_2(\alpha )$, and $T_2(\alpha )$ implies $T_1(\alpha )$, as in
Theorem 22.17 of [17].

In this section we consider probability measure spaces which describe
stochastic processes in discrete time. Let $\hat X =(X,d, \mu_1)$ be a
metric probability space and $\xi^{(n)}=(\xi_j)_{j=1}^n$ a stochastic process
with state space $X$, so there exists a probability $\mu_n$ on $X^n$
that specifies the joint distribution of $\xi^{(n)}$. We endow $X^n$
with the metric $d_n(x^{(n)}, y^{(n)})=(\sum_{j=1}^n d(x_j,
y_j)^2)^{1/2}$ for all $x^{(n)}=(x_j)_{j=1}^n$ and
$y^{(n)}=(y_j)_{j=1}^n$ in $X^n$. \par
\indent (ii) Whereas $LSI(\alpha )$ and $T_2(\alpha )$ imply spectral
gap inequalities, otherwise known as Poincar\'e inequalities, $T_1(\alpha
)$ inequalities do not in general; see [22, Theorem 22.17]. 
Hence Theorem 3.2 is less precise
than Proposition 3.1 and the results of [17].\par

 We write $x^{(j)}=(x_1, \dots ,x_j)\in X^j$ and let
$\mu_{j+1}(dx^{(j+1)})=p_{j+1}(dx_{j+1}\vert x^{(j)})\mu_j(dx^{(j)})$
be the disintegration of $\mu_{j+1}$ with respect to the measure on
$X^{j+1}$ that is pushed forward from $\mu_j$ by $\iota_j$. In other
words, $p_{j+1}(dx_{j+1})\vert x^{(j)})$ gives the distribution of
$\xi_{j+1}$ given $\xi^{(j)}=x^{(j)}$. \par

\indent Given a metric probability space that satisfies the transportation
inequality $T_1(\alpha )$, one can introduce potentials and
corresponding families of metric probability
spaces with similar properties. The following result is applicable to Gaussian
random Fourier series in subsequent sections.
\vskip.05in
\noindent {\bf Proposition 2.1} {\sl Suppose that $(\Omega , d,\mu )$ satisfies
$T_1(\alpha )$ for some $\alpha >0$ and that $V_n, V:\Omega\rightarrow [0,
\infty )$ are measurable potentials such that\par 
\indent (i) $Z(\lambda )=\int e^{\lambda V(u)}\mu (du)$ has $Z(\lambda
)<\infty$ for all $\lambda >0$;\par
\indent (ii) there exists $C>0$ such that $V_n(u)\leq CV(u)$ for all $u\in
\Omega$ and all $n$;\par
\indent (iii) $V_n(u)\rightarrow V(u)$ as $n\rightarrow\infty$ for all $u\in
\Omega$.\par
\noindent Then there exists $Z_n>0$ and $\varepsilon >0$ such that $d\nu_n =Z_n^{-1}e^{V_n}d\mu$ and
$d\nu =Z(1)^{-1}e^Vd\mu$ are probability measures on $\Omega$ that all satisfy
$T_1(\varepsilon )$ and $W_1(\nu_n, \nu )\rightarrow 0$ as
$n\rightarrow\infty.$}\par
\vskip.05in
\noindent {\bf Proof.} Let $Z_n=\int e^{V_n}d\mu$, so that $1\leq Z_n\leq
Z(C)$ and by the dominated convergence theorem $Z_n\rightarrow Z(1)$ as
$n\rightarrow\infty$. Hence $\nu_n\in {\hbox{Prob}}_0(\Omega )$. There exists
$u_0\in \Omega$ and $K(\mu )>0$ such that
$$\eqalignno{\int e^{\alpha^2 d(u,u_0)^2/2}\nu_n(du)&\leq \Bigl( \int Z_n^{-2}e^{2V_n(u)}
\mu (du)\Bigr)^{1/2}\Bigl( \int e^{\alpha^2 d(u,u_0)^2}\mu
(du)\Bigr)^{1/2}\cr
&\leq Z(2C)^{1/2}K(\mu )^{1/2},&(2.4)\cr}$$
\noindent for all $n$. Hence by Theorem 22.22, $\nu$ and all the $\nu_n$ satisfy
$T_1(\varepsilon )$, so in particular,
$$ W_1(\omega , \nu )\leq {{2}\over{\alpha}}\Bigl(1+(1/2)
\log Z(2C)K(\mu )\Bigr)^{1/2}{\hbox{Ent}}(\omega \mid \nu )^{1/2}\eqno(2.5)$$
\noindent for all $\omega\in {\hbox{Prob}}_0(\Omega )$ that are of finite
relative entropy with respect to $\nu$. Finally, we can take $\omega =\nu_n$
and note that by the dominated convergence theorem, ${\hbox{Ent}}(\nu_n\mid \nu )\rightarrow
0$ as $n\rightarrow\infty$, so $W_1(\nu_n, \nu )\rightarrow 0$.\par 
\vskip.05in
\vskip.1in
\noindent {\bf Proposition 2.1} {\sl (i) Suppose that there exists
$\alpha_0>0$ such that $\mu_1$ and $p_k(\, \cdot\, \vert x^{(k-1)})$
satisfy $T_2(\alpha_0)$ for all $x^{(k-1)}\in X^{(k-1)}$ and all $k=2,
\dots ,n$; }\par
\indent {\sl (ii) Suppose that the maps $x^{(k-1)}\mapsto p(\,\cdot\,
\vert x^{(k-1)})$ are Lipschitz from} $(X^{k-1}, d_{k-1})\rightarrow
({\hbox{Prob}}_2(X^{k-1}), d_{W_2})$ {\sl in the sense that there exist
$0\leq\rho_j<1$ and $0\leq R<1$ such that} 
$$d_{W_2}^2 (p_k(\, \cdot\, \vert x^{(k-1)}),p_k(\, \cdot\, \vert
y^{(k-1)})\leq\sum_{j=1}^{k-1} \rho_{k-j}d_j^2(x_j, y_j)\qquad (x_j,
y_j\in X),\eqno(2.6)$$
\noindent {\sl where $\sum_{j=1}^n \rho_j\leq R$.\par
\indent (iii) Suppose further that 
$$\sum_{j=1}^\infty \Bigl( \int_X {\hbox{var}}(X, d, p_{j+1}(\, \cdot
\, \vert x^{(j)})) \mu_j(dx^{(j)})\Bigr)^{1/2}<\infty .\eqno(2.7)$$
\indent (1) Then $(\hat X^j)$ converges in $D_{L^2}$ to
a metric probability space $\hat X^\infty$.}\par
\indent {(2) Also $\mu_n$ satisfies $T_2((1-\sqrt{R})^2\alpha_0)$, and
there exists $\varepsilon >0$ independent of $n$ such that}\par
$$\int\!\!\!\int_{X^n\times X^n} e^{\varepsilon d_n^2(x^{(n)},
y^{(n)})} \mu_n(dx^{(n)})\mu_n(dy^{(n)})<\infty .\eqno(2.8)$$
\vskip.05in
\noindent {\bf Proof.} (1)  We prove that the series
$\sum_{j=1}^\infty D_{L^2}(\hat X^j, \hat X^{j+1})$ converges and hence
$(\hat X^j)_{j=1}^\infty $ converges in the complete metric space $({\bf X}, D)$. For
$\varepsilon>0$, let $\theta_j:X^j\rightarrow X$ be a continuous
function such that
$${\hbox{var}}(X, d,p_{j+1}(\, \cdot \, \vert x^{(j)}))\leq \int_X
d^2(\theta_j(X^{(j)}), y_{j+1})p_{j+1}(dy_{j+1}\vert x^{(j)})
+2^{-j}\varepsilon.\eqno(2.9)$$
\indent We introduce a coupling $\delta_j$ of $(X^j, d_j)$ and $(X^{j+1}, d_{j+1})$ by 
$$\eqalignno{ \delta (x^{(j)}, y^{(j+1)})^2&=d_j^2(x^{(j)}, y^{(j)})+d^2(\theta_j(x^{(j)},
y_{j+1}),\qquad (x^{(j)}\in X^{j}, y^{(j+1)}\in X^{j+1});\cr
  \delta^2 (x^{(j)}, z^{(j)})&=d_j^2(x^{(j)}, z^{(j)})\qquad (x^{(j)}, y^{(j)}\in X^j);\cr
\delta^2 (y^{(j+1)},z^{(j+1)})&=d_j^2(y^{(j)},z^{(j)})+d^2(y_{j+1},z_{j+1})\qquad (y^{(j+1)}, z^{(j+1)}\in
X^{j+1}).&(2.10)\cr}$$ 
\indent We also observe that the probability measure on $X^j\times X^{j+1}$ given by
$$\pi_j(dx^{(j)}dy^{(j+1)})=p_{j+1}(dy_{j+1}\vert y^{(j)})\mu_j(dy^{(j)})
\delta_{y^{(j)}}(dx^{(j)})\eqno(2.11)$$ 
\noindent where $\delta$ here stands for the Dirac point mass, has marginals $\mu_j(dx^{(j)})$ and
$\mu_{j+1}(dy^{(j+1)}).$ Combining this with (), we deduce that
$$D^2(\hat X^j, \hat X^{j+1})\leq \int\!\!\!\int_{X^j\times X}d^2(\theta_j(x^{(j)}),
y_{j+1})p_{j+1}(dy_{j+1}\vert x^{(j)})\mu_j(dx^{(j)})+2^{-n}\varepsilon .\eqno(2.12)$$
\indent (2) By a particular case of Theorem 2.1 of
[4], all of the $\mu_n$ satisfy $T_2((1-\sqrt{R})^2\alpha_0)$. Further, 
$T_2((1-\sqrt{R})^2\alpha_0)$ implies $T_1((1-\sqrt{R})^2\alpha_0)$, which in turn is
equivalent to a Gaussian concentration inequality, as in (2.8). The details are described in 
[17, Theorem 22.22].\par

The proof mainly depends upon a concentration inequality for random
Fourier series. It is convenient to use scaled Fourier coefficients $x_j, y_j$
and write $z^{(n)}=(x_{-n}+iy_{-n}, \dots ,x_n+iy_n)$. Then we introduce the scalar potential
from the Hamiltonian $H_n$ by 
$$V(z^{(n)})=\int_{{\bf T}}\Bigl\vert\sum_{j=-n; j\neq 0}^n
{{z_je^{ij\theta}}\over{j}}\Bigr\vert^4 {{d\theta}\over{2\pi}}.\eqno(3.14)$$
\noindent We also introduce $\gamma^{(n)}(dz^{(n)})=\otimes_{j=-n}^n \gamma_j
(dx_j)\gamma'_j(dy_j)$ be the product Gaussian measure on ${\bf
C}^{2n+1}$. When $(z^{(n)})$ is randomly distributed according to $\gamma^{(n)}$, the potentials satisfy the following concentration inequality, with
constants that are independent of $n$.\par
\vskip.05in
\noindent {\bf Lemma 3.3} {\sl The following inequality holds, where the right-hand side is
bounded by constant independent of $n$,}
$$\eqalignno{\int_{\Omega_N} \exp&\Bigl({\lambda V(z^{(n)})}\Bigr)
 \gamma^{(n)}(dz^{(n)})\cr
&\leq
\sqrt{2}\exp\Bigl( 2^{16}3^{-3}\pi^2N\lambda^2\int
V(z^{(n)})\gamma^{(n)}(dz^{(n)})\Bigr)\qquad (\lambda\in {\bf R}).&(3.15)\cr}$$
\vskip.05in
\noindent {\bf Proof.} We observe that 
$$\eqalignno{ \Bigl\vert{{\partial V}\over{\partial
x_k}}\Bigr\vert&\leq{{4}\over{\vert k\vert}}\int_{{\bf T}} \Bigl\vert \sum_{j=-n; j\neq 0}^n
{{z_je^{ij\theta}}\over{j}}\Bigr\vert^3 {{d\theta}\over{2\pi}} \cr
&\leq {{4}\over{\vert k\vert}}\Bigl(\int_{{\bf T}} \Bigl\vert \sum_{j=-n; j\neq 0}^n
{{z_je^{ij\theta}}\over{j}}\Bigr\vert^4 {{d\theta}\over{2\pi}}\Bigr)^{1/2} 
\Bigl(\int_{{\bf T}} \Bigl\vert \sum_{j=-n; j\neq 0}^n
{{z_je^{ij\theta}}\over{j}}\Bigr\vert^2
{{d\theta}\over{2\pi}}\Bigr)^{1/2}&(3.16)\cr}$$
\noindent where the final factor is bounded by $\sqrt{N}$ for fields over $\Omega_N$;
hence the gradient $\nabla V=(\partial V/\partial x_k; \partial V/\partial
y_k)_{k=-\infty}^\infty$ satisfies $\Vert \nabla V\Vert_{\ell^2}\leq 4\pi
V(z^{(n)})\sqrt{(N/3)},$ for all $z^{n}$ corresponding to fields over $\Omega_N$. We deduce that $V(z^{(n)})^{1/2}$ is $4\pi
\sqrt{(N/3)}$-Lipschitz with respect to the standard $\ell^2$ metric, and hence by the
concentration inequality for the Gaussian measure 
$$\int_{\Omega_N} \exp\Bigl({t\sqrt{\vert\lambda\vert} V(z^{(n)})^{1/2}}\Bigr) \gamma^{(n)}
(dz^{(n)})$$
$$\leq
\sqrt{2}\exp\Bigl( 8\pi^2\vert\lambda\vert
t^2/3+t\sqrt{\vert\lambda\vert}\int
V(z^{(n)})\gamma^{(n)}(dz^{(n)})\Bigr)\qquad (t>0).\eqno(3.17)$$
\noindent We select $\sigma^2=3/(32\pi^2N\vert\lambda\vert )$, and integrate the
preceding inequality with respect to $e^{-t^2/2\sigma^2}dt/\sqrt{2\pi\sigma^2}$, thus
obtaining
$$\eqalignno{\int_{\Omega_N} e^{\sigma^2\vert\lambda\vert
V(z^{(n)})/2}\gamma^{(n)}(dz^{(n)})&= \int_{\Omega_N} 
\int_{-\infty}^\infty e^{t\sqrt\vert\lambda\vert V(z^{(n)})^{1/2}/2 -t^2/2\sigma^2}
{{dt}\over{\sqrt{2\pi\sigma^2}}} \gamma^{(n)}(z^{(n)})\cr
&\leq \sqrt{2} \exp\Bigl( 16\lambda^2\pi^2N\int_{\Omega_N}
V(z^{(n)})\gamma^{(n)}(dz^{(n)})/3\Bigr).&(3.18)\cr}$$
\noindent By the M. Riesz theorem, there exists a constant $C$ such that $V(z^{(n)})\leq
CV(z^{(\infty )})$ whenever $z^{(\infty)}=(z_j)_{j=1}^\infty $ is a bounded complex sequence that has 
initial segment $z^{(n)}=(z_j)_{j=1}^n$. Hence 
$$\int V(z^{(n)})\gamma^{(n)}(dz^{(n)})\leq C\int V(z^{(\infty )})\gamma^{(\infty )}(dz^{(\infty
)}),\eqno(3.19)$$
\noindent where the right-hand side is finite by basic results about Brownian loop.\par
By a rescaling of $\lambda$, we obtain the inequality of the Lemma.\par
\rightline{$\square$}\par
\vskip.05in
\noindent {\bf Proof of Theorem 3.2.} On $X^n=P_nL^2\cap \Omega_N$, the Gibbs  measure is expressed in
Fourier coefficients as 
$$\eqalignno{\nu_n (da^{(n)}db^{(n)})&=Z(\lambda ,N,n)^{-1}{\bf I}_{\{\sum_{j=-n}^n
(a_j^2+b_j^2)\leq N\}}(a^{(n)};b^{(n)})\exp\Bigl( \lambda\int_{\bf
T}\bigl\vert\sum_{j=-n}^n
(a_j+ib_j)e^{ij\theta}\bigr\vert^4{{d\theta}\over{2\pi}}\Bigr)\cr
&\prod_{j=-n;j\neq 0}^n
j^2e^{-j^2(a_j^2+b_j^2)/2}{{da_jdb_j}\over{2\pi}},&(3.20)\cr}$$
\noindent where $Z(N,\lambda,n)>0$ is a normalizing constant. By Lemma 3.3, there exist
constants $\kappa (N, \lambda )$ and $K(N, \lambda )$ independent of $n$ such that
$$0< \kappa (N, \lambda )\leq Z(N, \lambda ,n)\leq K(N, \lambda )\qquad (N>0, \beta\in {\bf
R}).\eqno(3.21)$$
\indent We create a  family of maps 
$X^n\rightarrow X^{n+\ell}$ by $(a^{(n)}; b^{(n)})\mapsto (a^{(n+\ell)}; b^{(n+\ell )}),$ where
$$(a^{(n)}) \mapsto (a_{-n-\ell},\dots a_{-n-1},a^{(n)}, 
a_{n+1}, \dots a_{n+\ell})\eqno(3.22)$$
\noindent depending upon the new variables
$$a_{-n-\ell},\dots , a_{-n-1},a_{n+1},\dots
,a_{n+\ell};b_{-n-\ell},\dots , b_{-n-1},b_{n+1}, \dots b_{n+\ell}$$
\noindent which are jointly distributed according to
the transition kernel $q_{n+1, n+\ell}$ that is uniquely specified by 
$$\nu_{n+\ell}(da^{(n+\ell)}db^{(n+\ell)})=dq_{n+1,
n+\ell}(\,\cdot\, \,\mid a^{(n)};b^{(n)})\nu_{n}(da^{(n)}db^{(n)}).\eqno(3.23)$$
\noindent With the Gaussian density 
$$g_{n+1, n+\ell}=\prod_{j=n+1}^{n+\ell}j^4e^{-j^2(a_{-j}^2+a_{j}^2+
b_{-j}^2+
b_{j}^2)/2}/(2\pi )^2,\eqno(3.24)$$
\noindent on ${\bf R}^{4\ell}$, Sturm [20, Corollary 4.21] has shown that the sequence of metric measure spaces 
$({\bf R}^{4\ell}, \ell^2, g_{n+1,n+\ell})$
converges in $D_{L^2}$ as $\ell\rightarrow\infty$ to $({\bf R}^\infty; \ell^2, \Gamma_n )$, where
$\Gamma_n$ is the Gaussian measure with auto-covariance matrix ${\hbox{diag}}(I_4/(n+1)^2,
I_4/(n+2)^2,\dots )$. 
 We have a Cauchy--Schwarz inequality
$$\eqalignno{ \int_{{\bf
R}^{4\ell}}&\Bigl(\sum_{j=n+1}^{n+\ell}(a_{-j}^2+a_{j}^2+b_{-j}^2+b_{j}^2)\Bigr)
dq_{n+1,n+\ell}(\,\,\mid a^{(n)};b^{(n)})&(3.25)\cr
&\leq\Bigl(\int_{{\bf
R}^{4\ell}}\Bigl(
\sum_{j=n+1}^{n+\ell}(a_{-j}^2+a_{j}^2+b_{-j}^2+b_{j}^2)\Bigr)^2
dg_{n+1, n+\ell}\Bigr)^{1/2}\Bigl( \int_{{\bf R}^{4\ell}}\Bigl(
{{dq_{n+1, n+\ell}}\over{dg_{n+1, n+\ell}}}\Bigr)^2dg_{n+1,
n+\ell}\Bigr)^{1/2}.\cr}$$
\noindent When we multiply out the first squared sum, we obtain several expressions which resemble
$$\eqalignno{\int_{{\bf R}^{4\ell}}\Bigl( \sum_{j,k=n+1; j\neq k}^{n+\ell} a_j^2a_k^2+\sum_{j=n+1}^{n+\ell}
a_j^4\Bigr) dg_{n+1,n+\ell}&=\sum_{j,k=n+1; j\neq k}^{n+\ell}
{{1}\over{j^2k^2}}+\sum_{j=n+1}^{n+\ell}{{3}\over{j^4}}\cr
&\leq 3\Bigl(\sum_{j=n+1}^{n+\ell} {{1}\over{j^2}}\Bigr)^2\leq {{C}\over{n^2}}.&(3.26)\cr}$$
\noindent Integrating with respect to $\nu_n$, and applying Cauchy--Schwarz
again, we obtain
$$\eqalignno{&\int_{{\bf
R}^{4(n+\ell )}}\Bigl(\sum_{j=n+1}^{n+\ell} (a_{-j}^2+a_{j}^2+b_{-j}^2+b_{j}^2)\Bigr)\nu_{n+\ell}
(da^{(n+\ell)}db^{(n+\ell)})\cr
&\leq {{C}\over{n+1}}\Bigl( 
\int_{{\bf R}^{4(n+\ell )}}\Bigl(
{{dq_{n+1, n+\ell}}\over{dg_{n+1, n+\ell }}}\Bigr)^2dg_{n+1,
n+\ell}\nu_n(da^{(n)}db^{(n)})\Bigr)^{1/2}\cr
&\leq {{C}\over{n+1}}{{Z(N,\lambda .n)}\over{Z(N, \lambda ,n+\ell )}}\Bigl( 
\int_{{\bf R}^{4(n+\ell )}}\exp\Bigl( 2\lambda V(z^{(n+\ell )})-2\lambda V(z^{(n)})\Bigr)\nu_{n+\ell
}(dz^{(n+\ell )})\Bigr)^{1/2}.&(3.27)\cr}$$
\indent From Lemma 3.3 and (3.21), we deduce that there exists $C(N, \lambda )>0$ such that
$${\hbox{D}}_{L^2}(\hat X^{n+\ell}, \hat X^n)\leq {{C(N, \lambda )}\over{n+1}}\qquad 
(n,\ell\geq 1),\eqno(3.28)$$
\noindent so that $(\hat X)_{n=1}^\infty$ gives a 
Cauchy sequence in $({\bf X}_1, {\hbox{D}}_{L^2})$,
converging to $(\Omega_N, \Vert\, \cdot\, \Vert_{L^2}, \nu_N^\lambda )$.\par
\rightline{$\square$}\par
\vskip.05in

\indent Suppose that $u$ is a solution of $KdV$ with initial value $u(x,0)$ in
the support of $\nu_N^\lambda$. The Hill's operator
$$L_u=-{{d^2}\over{dx^2}}+u\eqno(4.9)$$
\noindent is densely defined and self-adjoint in $L^2({\bf R}; {\bf C})$ with
periodic spectrum $(\lambda_j)_{j=1}^\infty $, where
$\lambda_0<\lambda_1\leq\lambda_2<\dots $ satisfies
$\lambda_j\rightarrow\infty$ as $j\rightarrow\infty$. Each $\lambda_j$ does not
change with $t$ as $u$ undergoes the flow associated with $KdV$.\par
\indent Let $P_n:L^2\rightarrow {\hbox{span}}\{ \cos jx, \sin jx: j=0, \dots
,n\}$ be the Dirichlet projection and 
$$H_n(u)={{1}\over{2}}\int_{\bf T} \Bigl( D_n{{\partial u}\over{\partial x}}
(x,t)\Bigr)^2 {{dx}\over{2\pi}} -{{\lambda}\over{6}}\int_{\bf T} (D_nu(x,t))^3
{{dx}\over{2\pi}},\eqno(4.10)$$
where $\lambda>0$ is the reciprocal temperature. Then the canonical equation of
motion ${{\partial u}\over{\partial t}}=
{{\partial }\over{\partial x}}{{\delta H}\over{\delta u}}$ gives the KdV
equation
$${{\partial u}\over{\partial t}}=-{{\partial^3 u}\over{\partial
x^3}}-\lambda D_n(u{{\partial u}\over{\partial x}}).\eqno(4.11)$$
Let $L_{u_n}=-{{d^2}\over{dx^2}}+u_n(x,t)$ be the corresponding family of
Hill's operators in $L^2({\bf T})$. The periodic spectrum of $L_{u_n}$ is
$(\lambda_j(u_n (\, ,t)))_{j=0}^\infty$; note that $\lambda_j(u_n(\, t))$
varies as $u_n$ undergoes the evolution generated by the Hamiltonian flow. \par  

\vskip.05in
\noindent {\bf Corollary 4.3} {\sl For each $j$, the periodic eigenvalue   
$\lambda_j(u_n (\, ,t))$ gives a random variable on $\hat X^n$, such that 
$\lambda_j(u_n (\, ,t))\rightarrow \lambda_j(u)$ as $n\rightarrow\infty$.}
\vskip.05in
\noindent {\bf Proof.} The eigenvalues are continuous with respect to the
potential in the $L^2$ norm topology.\par

\noindent {\bf Lemma 6.1} {\sl Let the metric probability space $(X, d, \mu )$ satisfy $T_1(\alpha )$ and let
$F:X\rightarrow {\bf R}$ be $L$-Lipschitz. Then}
$$\int_Xe^{\alpha F(\xi )^2/(4L^2)} \mu (d\xi )\leq \sqrt{2}e^{\alpha (\int Fd\mu
)^2/(2L^2)}.\eqno(6.7)$$
\vskip.05in
\noindent {\bf Proof.} The transportation inequality implies the Gaussian concentration inequality 
$$\int_X e^{tF(\xi )}\mu (d\xi)\leq e^{t\int Fd\mu +L^2t^2/(2\alpha )}
\qquad (t\in {\bf R}),\eqno(6.8)$$
\noindent so we choose $\sigma^2=\alpha /(2L^2)$ and integrate with respect to 
the Gaussian density\par
\noindent $e^{-t^2/(2\sigma^2)} /\sqrt{2\pi\sigma^2}$. The result emerges once we change the
order of integration.\par
\vskip.05in
\noindent {\bf Proposition 6.2} {\sl There exist constants $\kappa_1, \kappa_2, \kappa_3>0$
such that for all $N_n\leq \kappa_1\log n$ and
$0<\lambda_n<\kappa_2/(\log n)^2$, the metric
probability space $(X^n, L^2, \nu_n)$ satisfies $T_1(\kappa_3/N_n).$}\par
\vskip.05in
\noindent {\bf Proof.} By (22.16) of Villani, all $\omega_n\in {\hbox{Prob}}_1(X^n)$ that are
of finite relative entropy with respect to $\nu_n$ satisfy
$$W_1(\omega_n, \nu_n)\leq \sqrt{ {{2}\over{\varepsilon_n}}{\hbox{Ent}}(\omega_n\mid
\nu_n)}
\Bigl( 1+\log \int_{\Omega_n} e^{\varepsilon_n\Vert
u\Vert^2_{L^2}}\nu_n(du)\Bigr)^{1/2},\eqno(6.9)$$
\noindent  for all $\varepsilon_n>0$ such that the final integral converges. So it suffices to
show that selecting $0<\varepsilon_n<\kappa_3/N_n$, the factor in parenthesis is bounded
independent of $n$. \par
\indent In the Lemma 6.1, we choose $X=\Omega_n$ and 
$$\mu_n (dxdy)=\zeta_n^{-1}\prod__{k=(k_1,k_2)\in {\bf Z}^2; \vert k_1\vert, \vert k_2\vert \leq n} e^{-(x_k^2+y_k^2)/2}
{{dx_kdy_k}\over{2\pi}}\eqno(6.10)$$
\noindent where the normalizing constant satisfies $\zeta_n\rightarrow 1$ as
$n\rightarrow\infty$, and $\alpha_n\rightarrow\infty$ by standard results about Gaussian
measures. Then we choose $F:X^n\rightarrow {\bf R}$ as
$$F(x, y)=\Bigl( \int_{{\bf T}^2} \Bigl\vert \sum__{k=(k_1,k_2)\in {\bf Z}^2; \vert k_1\vert, \vert k_2\vert \leq n}
{{(x_k+iy_k)e^{ik\cdot \theta}}\over{\vert k\vert}}\Bigr\vert^4 {{d^2\theta}\over{(2\pi
)^2}}\Bigr)^{1/2},\eqno(6.11)$$
\noindent where $x=(x_k)$ and $y=(y_k)$ are such that $u\in \Omega_n$. Then by repeating the
computation of section 4, we see that $F$ is $L$-Lipschitz where $L\leq N_n$. \par

\indent  A natural scaling is $N_n=C\log n$ and $\lambda_n=c/(\log n)^2$; indeed, letting
$(\gamma_k, \gamma_k')_{k\in {\bf Z}^2}$ be standard mutually independent Gaussian random
variables, we note that the expectation of 
$$ \int_{{\bf T}^2}\Bigl\vert \sum__{k=(k_1,k_2)\in {\bf Z}^2; \vert k_1\vert, \vert k_2\vert \leq n} {{e^{ik\cdot\theta}
(\gamma_k+i\gamma_k')}\over{\vert k\vert}}\Bigr\vert^2   
{{d\theta_1}\over{2\pi}} {{d\theta_1}\over{2\pi}}\eqno(6.12)$$
\noindent is of order $\log n$ as $n\rightarrow\infty$.\par 
\indent Now the principal measure on $X^n$ is
$d\nu_n=Z_n(\lambda_n)^{-1}e^{\lambda_n\int \vert
u\vert^4} d\mu_n$, where the normalizing constant satisfies 
$$Z_n(\lambda )=\int_{\Omega_n} e^{\lambda\int \vert u\vert^4}\mu_n(du)\eqno(6.13)$$
\noindent so $1\leq Z_n(\lambda_n)$ by Jensen's
inequality and by () we can choose $\kappa_2, C>0$ such that 
$Z_n(\lambda )\leq C$ for all $0<\lambda_n<2\kappa_2/(\log n)^2.$\par 
\indent Finally, we compare $\nu_n$ and $\mu_n$ by using the Cauchy--Schwarz inequality
$$\eqalignno{ \int_{\Omega_n}& e^{\varepsilon_n\Vert u\Vert^2_{L^2}}\nu_n(du)\cr
&\leq \Bigl(  \int_{\Omega_n} e^{2\varepsilon_n\Vert
u\Vert^2_{L^2}}\mu_n(du)\Bigr)^{1/2}Z_n(\lambda_n)^{-1} 
\Bigl(  \int_{\Omega_n}e^{2\lambda_n\int\vert
u\Vert^4}\mu_n(du)\Bigr)^{1/2}.&(6.14)\cr}$$
\noindent We recognise the final integral as $Z_n(2\lambda_n)\leq C$, while we can control the
first integral by taking $\varepsilon_n=1/N_n$.\par\vskip.1in

\indent {\bf 6. Logarithmic Sobolev inequalities for modified ensembles in 1D}

\indent Lebowitz Rose and Speer [13] consider possible modifications to the interaction energy $\int
\vert \phi\vert^pdx$ which prevent the field from collapsing. For the $KdV$ equation and $NLS$,
we consider replacing $(\lambda /k)\int u(x)^k dx$ by 
$${{\lambda }\over{k}}\int_{\bf T} u(x)^k\, dx-{{\lambda^2}\over{2}}\int_{\bf T} u(x)^{2(k-1)}
\d, dx.$$
\noindent With this modified interaction terms, the corresponding Hamiltonian generates some
natural partial differential equations. The corresponding Gibbs measure can be normalized,
since the defocussing term involving $\lambda^2$ stabilizes the focussing term in $\lambda$;
indeed, we show that the Gibbs measure for $k=3$ satisfies a logarithmic Sobolev inequality.\par
\indent For integers $k=2,4, \dots $, we introduce the modified ensemble which has Hamiltonian
$$ H(u)={{1}\over{2}}\int_{\bf T}\Bigl({{\partial u}\over{\partial x}}\Bigr)^2
dx-{{\lambda}\over{k}}\int_{\bf T} u(x)^k dx+{{\lambda^2}\over{2}}
\int_{\bf T} u(x)^{2(k-1)}dx.\eqno(6.1)$$
The canonical equation 
$${{\partial u}\over{\partial t}}={{\partial }\over{\partial x}}{{\delta H}\over{\delta u}}$$
gives the differential equation
$${{\partial u}\over{\partial t}}=-{{\partial^3u}\over{\partial x^3}}-
\lambda (k-1)u^{k-2}{{\partial
u}\over{\partial x}}+{{\lambda^2}\over{2}}(2k-3)u^{2k-4}
{{\partial u}\over{\partial x}},\eqno(6.2)$$
so that for suitably differentiable real periodic functions $u$, 
both $\int_{\bf T} u(x)^2dx$ and $H(u)$ are
invariant with respect to time.\par
\indent So we introduce the formal Gibbs measure
$$Z^{-1}{\bf I}_{\Omega_N}(u) e^{-H(u)} \prod_{x\in [0, 2\pi ]} du(x).\eqno(6.3)$$
\vskip.05in
\noindent {\bf Proposition 6.1} {\sl For $k=2,4$ and all $N$, the Gibbs measure can be normalized to give a Radon probability
measure on $\Omega_N$. This Gibbs measure satisfies a logarithmic Sobolev
inequality.}\par
\vskip.05in
\noindent {\bf Proof.} Let $\mu$ be the probability measure on $\Omega_N$ that is induced by the map
$\omega\rightarrow \sum_{k=1}^\infty (\gamma_k\cos k\theta +\gamma'_k\sin k\theta )/k$
where $(\gamma_k, \gamma_k')_{k=1}^\infty$ are mutually independent standard Gaussian random variables.\par 

\indent  Let $V:L^2\rightarrow [0,\infty ]$ be the convex function
$$V(u)={{1}\over{k}}\int_{\bf T} u(x)^k dx,\eqno(6.4)$$
\noindent and let $c:L^2\times L^2\rightarrow [0, \infty )$ be the cost
function $c(u,w)=2^{-1}\Vert u-w\Vert^2_{L^2}.$ Then we introduce the Hopf--Lax functional
$$QV(w)=\inf\{ V(u)+c(u,w): u\in L^2\}.\eqno(6.5)$$

\noindent Then $\mu$ satisfies a logarithmic Sobolev inequality, so by a theorem of Bobkov,
Gentil and Ledoux [6, Corollary 6.1] 
$$\int_{\Omega_N} e^{QV(w)}\mu (dw)\leq e^{\int V(w)\mu (dw)}.\eqno(6.6)$$
As in Theorem 4.1, the measure $\mu$ is a limit of measures on the finite-dimensional
spaces $X_n$. Now when $w$ is a trigonometric polynomial, 
$V$ is differentiable at $w$, and we can form $\nabla V$ by partially
differentiating with respect to the Fourier coefficients. Hence by convexity, we have
$$V(u)\geq V(w)+\langle \nabla V(w), u-w\rangle_{L^2},\eqno(6.7)$$
so that 
$$QV(w)\geq V(w)-2^{-1}\Vert \nabla V(w)\Vert^2.\eqno(6.8)$$
Hence we have the inequality
$$\int_{\Omega_N} \exp\Bigl(    
{{\lambda}\over{k}}\int_{\bf T} u(x)^k
dx-{{\lambda^2}\over{2}}\int_{\bf T} 
u(x)^{2(k-1)}dx\Bigr)\mu
(du)\leq \exp \int_{\Omega_N}\Bigl({{\lambda}\over{k}}\int_{\bf T} u(x)^k dx\Bigr)
\mu (du),\eqno(6.9)$$
where the right-hand side is finite by elementary results about random Fourier series.
\indent We consider the functional
$$W(u)={{\lambda^2}\over{2}}\int_{\bf T}u(x)^4\, dx-{{\lambda}\over{3}}\int_{\bf
T} u(x)^3 \, dx,\eqno(6.10)$$
\noindent and show that $W(u)+5\Vert u\Vert^2_{L^2}/9$ is uniformly convex.
Indeed, for $0<t<1$, we have
$$\eqalignno{{{\lambda^2}\over{2}}&\Bigl( tu^4+(1-t)v^4-\bigl( tu+(1-t)v\bigr)^4\Bigr)
-{{\lambda}\over{3}}\Bigl( tu^3+(1-t)v^3-\bigl(tu+(1-t)v\bigr)\Bigr)\cr
&=t(1-t)(u-v)^2\Bigl({{\lambda^2}\over{2}}\bigl((1+t+t^2)u^2
+2(1+t-t^2)uv+(3-3t+t^2)v^2\bigr)\cr
&\qquad  -{{\lambda}\over{3}}\bigl(
(1+t)u+(2-t)v\bigr)\Bigr).&(6.11)\cr}$$
To analyse the quadratic form in this, we introduce the matrix
$$A=\left[\matrix{t^2+t+1& -t^2+t+1\cr -t^2+t+1& t^2-3t+3\cr}\right]\eqno(6.12)$$
\noindent which is positive definite with eigenvalues
$$\lambda_{\pm}=t^2-t+1\pm\sqrt{ (t^2+1)(t^2-2t+2)}.\eqno(6.13)$$
\noindent Hence $\det A$ attains its smallest value $3/2$ at $t=1/2$, and
$\lambda_+$ attains its smallest value $3$ at $t=1/2$; also $\lambda_-$ attains
its smallest value $1/2$ there. Hence 
$$\eqalignno{ t(1-t)(u-v)^2&\Bigl({{\lambda^2}\over{4}}\bigl(u^2+v^2\bigr)
 -{{\lambda}\over{3}}\bigl((1+t)u+(2-t)v\bigr)\Bigr)&(6.14)\cr
&= t(1-t)(u-v)^2\Bigl( \bigl( {{\lambda u}\over{2}}-{{1+t}\over{3}}\bigr)^2+
\bigl( {{\lambda v}\over{2}}-{{2-t}\over{3}}\bigr)^2-
{{1}\over{9}}\bigl( (1+t)^2+(2-t)^2\bigr)\Bigr)\cr}$$
We deduce that 
$$\eqalignno{tH(u)+(1-t)H(v)&-H(tu+(1-t)v)\cr
&\geq {{t(1-t)}\over{2}}\int_{\bf T}\Bigl(
{{\partial u}\over{\partial x}}-{{\partial v}\over{\partial x}}\Bigr)^2
dx-{{5t(1-t)}\over{9}}\int_{\bf T} (u-v)^2dx.&(6.15)\cr}$$ 
This suffices to show that the Gibbs measure satisfies a logarithmic Sobolev
inequality, since $\Vert u\Vert^2_{L^2}$ is bounded on $\Omega_N$.\par
\vskip.05in
 A
more drastic change would involve the integrand $\lambda u^k/k-\lambda^2 u^{2(k-1)}
-(k-1)u^2/(4k-6)$, which is concave, and hence has a unique maximum at $u=0$.\par  
1}$ 

\noindent where the perturbing potential satisfies
$$\int_{\Omega_N}\exp\Bigl(\kappa \Vert\nabla
W_p(u)\Vert^2_{L^2}\Bigr)\mu_N(du)<\infty\eqno(7.12)$$     
\noindent for some  $\kappa>0$; 

indeed, for $2(p-1)<6$, Lebowitz, Rose and Speer have shown that the 
Gibbs measure may be normalized for $H_{2(p-1)}$ on $\Omega_N$, so by the inequality (2.16), the integral (2.18) 
is
finite for all $\kappa$. Hence $\nu_N$ satisfies $LSI(\alpha )$ for some 
$\alpha >0$ by Aida and Shigekawa's criterion [1, p. 154] and [17, Remark 21.5].\par 
\indent (ii) For $p=4$, [] proved that there exists $N_0>0$ such that the Gibbs measure 
associated with $H_{2(p-1)}$ can be normalized on $\Omega_N$ for all $N<N_0$, so (2.18) is
finite for some $\kappa>0$. Hence by choosing $\lambda >0$ sufficiently small, we can still satisfy the criterion of [1] and obtain
$LSI(\alpha)$, for some $\alpha>0$. Aida and Shigekawa deduce their logarithmic Sobolev inequality
[1, Theorem 1.1] from a spectral gap inequality and a defective logarithmic Sobolev
inequality which involves additional constants and a subtle
perturbation argument. Hence their result does
not appear to be a route towards a $LSI(\alpha )$ with sharp
constants in critical cases.\par \par
\indent (iii) Let $\mu$ be as in (i), so that $\mu$ satisfies $LSI(\alpha_0 )$;

\noindent {\bf Proposition 8.1} {\sl For $0<s<1/2$ and $K>0$, let 
$$\Omega_{\infty
,K}=\{ u\in H^{-s}({\bf T}^2: {\bf C}): \sum_{k\in {\bf Z}^2\setminus \{
0\}}:\vert k\vert^{-2s}\vert \hat u(k)\vert^2\leq K\}.\eqno(8.9)$$
 Then there exist
$\kappa>0$ and $C>0$ such that} 
$$\int_{\Omega_{\infty ,K}}\exp\Bigl(\kappa \sum_{r=1}^\infty d_r^2\Bigr)\mu
(du)\leq C.\eqno(8.10)$$
\vskip.05in
\noindent {\bf Proof.} For each $r$, we use the index set $\{j\in {\bf
Z}^2\setminus \{ 0,-m\}; r-1<\vert j\vert \leq
r\}$. The map $(\gamma_j)\mapsto d_r$ is Lipschitz with respect
to the $\ell^2$ norm on those $(\gamma_j)$ such that
$\sum_j\gamma_je^{ij\cdot \theta}/\vert j\vert\in \Omega_{\infty ,K}$, so that 
$$\eqalignno{\vert d_r((\gamma_j)_j)&-d_r((\gamma_j')_j)\vert\cr
&\leq 
\Bigl( \sum_j {{\vert\gamma_j-\gamma_j'\vert^2}\over
{\vert j\vert^2}}\Bigr)^{1/2}\Bigl( \sum_j{{\vert\gamma_{j+m}\vert^2}
\over{\vert j+m\vert^2}}\Bigr)^{1/2}+\Bigl( \sum_j{{\gamma_j'\vert^2}\over{\vert j\vert^2}}\Bigr)^{1/2}   
\Bigl( \sum_j{{\vert\gamma_{j+m}-\gamma'_{j+m}\vert^2}\over{\vert
j+m\vert^2}}\Bigr)^{1/2}\cr
&\leq {{\sqrt{K}(r-\vert m\vert)^s}\over{r}}\Bigl(   
\sum_j\vert\gamma_j-\gamma_j'\vert^2\Bigr)^{1/2} +
{{\sqrt{K}(r^s}\over{r-\vert m\vert}}\Bigl(   
\sum_j \vert\gamma_{j+m}-\gamma_{j+m}'\vert^2\Bigr)^{1/2}.&(8.11)\cr}$$
\noindent Now the Gaussian measure $\mu$ satisfies $LSI(1)$, and hence the
concentration inequality
$$\int_{\Omega_{\infty , K}}\exp \bigl( \lambda d_r\bigr) d\mu \leq \mu
(\Omega_{\infty , K}) \exp\bigl({c_1\lambda^2K/r^{2-2s}}\bigr)\qquad
(\lambda\in {\bf R})\eqno(8.12)$$
\noindent for some $c_1>0$; next we integrate with respect to the Gaussian with
mean zero and variance $\sigma^2< r^{2-2s}/2c_1K$, obtaining
$$\int_{\Omega_{\infty , K}}\exp \Bigl( \sigma_r^2 d^2_r/2\Bigr) d\mu\leq 
\mu
(\Omega_{\infty , K})\Bigl(
1-{{2c_1K\sigma_r^2}\over{r^{2-2s}}}\Bigr)^{-1/2}.\eqno(8.13)$$
\noindent Finally, we use H\"older's inequality with $\sum_r 1/\sigma_r^2=1$ to
obtain
$$\int_{\Omega_{\infty ,K}}\exp\Bigl(\kappa \sum_{r=1}^\infty d_r^2\Bigr)\mu
(du)\leq\prod_{r=1}^\infty 
\Bigl(\int_{\Omega_{\infty ,K}}\exp\Bigl(\kappa \sigma_r^2 d_r^2\Bigr)\mu
(du)\Bigr)^{1/\sigma_r^2},\eqno(8.14)$$
\noindent where the right-hand side is convergent by ().\par

Under the Gaussian measure $\mu$, we have a martingale
difference sequence $(d_r)_{r=1}^\infty$, such that $\sum_{r=1}^\infty
d_r$ is convergent in $L^2$. Whereas the $d_j$ are not bounded,
nevertheless, Tao and
Vu [19] have shown that such martingales satisfy a concentration inequality
$$\mu\Bigl[ \bigl\vert\sum_{r=\ell}^Ld_r\bigr\vert
\geq\bigl(\sum_{r=\ell}^L \delta_r^2\bigr)^{1/2}\xi\Bigr] \leq
c_1e^{-c_2\xi^2}+\sum_{r=\ell}^L\mu \bigl[ \vert d_r\vert\geq
\delta_r\bigr]\eqno(8.15)$$
\noindent for some $c_1, c_2>0$. For $r\geq \ell \geq \vert m\vert+1$,
we choose $\delta_r=1/r^{\kappa +1/2}$ where $1/2-s>\kappa >0$, and observe 
that
$$\mu \Bigl([\vert d_r\vert>\delta_r]\cap \Omega_{\infty , K}\Bigr)\leq 
2e^{-4r^{1-2s-2\kappa}},\eqno(8.16)$$    

  Furthermore, by () we
have
$$\eqalignno{\Bigl(\int_{\Omega} \bigl\vert \widehat{(\vert u\vert)^2}
(m)-
\widehat{(\vert P_nu\vert^2)}(m)\bigr\vert^2\mu
(du)\Bigr)^{1/2}&=\Bigl(\int_\Omega \Bigl\vert
\sum_{r=n+1}^\infty d_r^{(m)}\Bigr\vert^2\mu (du)\Bigr)^{1/2}\cr
&\leq \Bigl({\bf E} \Bigl\vert\sum_{r=n+1}^\infty
\varepsilon_r\delta_r\Bigr\vert^2\Bigr)^{1/2}\cr
&=\Bigl(\sum_{r=n+1}^\infty \delta_r^2\Bigr)^{1/2}\leq {{C_0K_1^2}\over {n^\varepsilon\sqrt{\varepsilon}}},&(8.39)\cr}$$
\noindent hence $\widehat{(\vert P_nu\vert^2)}(m)$ converges to 
$\widehat{(\vert u\vert^2)}(m)$ almost surely, in $L^2$ and hence in $\mu$ 
measure as $n\rightarrow\infty$ for all $m$; essentially this follows from () by
the $L^2$-martingale convergence theorem.\par
\indent By standard results, convergence in $\mu$-measure and uniform
integrability together imply convergence in $L^1(\mu )$, each summand in ()
converges to zero in measure. Uniform integrability
follows from the inequality 
$$\eqalignno{\int_\Omega \exp\Bigl[ \sum_{m: \vert m\vert \geq r}&\Re  
\hat V(m)\Bigl(\widehat{\bigl( \vert P_nu\vert^2\bigr)}(m)-\widehat{\bigl(\vert
u\vert^2\bigr)}(m)\Bigr)\widehat{\bigl( \vert
u\vert^2\bigr)}(-m)\Bigr]\mu (du)\cr
& \leq \Bigl( \int_\Omega \exp\Bigl[ 2\lambda K_3\sum_{m: \vert m\vert \geq r}
{{\bigl\vert \widehat{(\vert u\vert^2)}(m)-\widehat{(\vert
P_nu\vert^2)}(m)\bigr\vert^2}\over{\vert
m\vert^{2+2s}}}\Bigr] \mu (du)\Bigr)^{1/2}\cr
&\quad\times \Bigl( \int_\Omega \exp\Bigl[ 2\lambda K_3\sum_{m: \vert m\vert \geq r}
{{\bigl\vert \widehat{(\vert u\vert^2)}(-m)\bigr\vert^2}\over{\vert
m\vert^{2+2s}}}\Bigr] \mu (du)\Bigr)^{1/2}&(8.40)\cr}$$
\noindent By Lemma 8.3, for all $\varepsilon >0$, there exists $r_0(\varepsilon )$ such that
all the integrals in () are less than $1+\varepsilon $ for all $n$ and all
$r\geq r_0(\varepsilon )$. Hence $W(P_nu )$ converges to $W(u)$ in $L^1(\mu
;\Omega)$ as $n\rightarrow\infty.$\par

The variance of $\hat X$ is
$${\hbox{var}}(\hat X)=\inf_z\Bigl\{
\int_X\delta_2(z,y)^2\mu_2(dy)\Bigr\}\eqno(2.5)$$
\noindent where the infimum is taken over all $\hat Y=(Y,\delta_2, \mu_2)$
that are isomorphic to $\hat X$, and all $z\in Y$.\par
\vskip.05in
\indent Evidently,
${\hbox{var}}(\hat X)$ is finite whenever $\mu_1\in
{\hbox{Prob}}_2(X)$.